\documentclass[11pt,fleqn]{article}
\RequirePackage[OT1]{fontenc}
\RequirePackage{amsthm,amsmath}
\RequirePackage[round]{natbib}
\RequirePackage[colorlinks,citecolor=blue,urlcolor=blue]{hyperref}

\usepackage{amsmath,amsfonts,amssymb,euscript,mathrsfs}
\usepackage{slashbox}
\usepackage{float}
\usepackage[final]{pdfpages}

\usepackage{txfonts}
\usepackage[T1]{fontenc} 
\usepackage{fullpage}
\usepackage{graphicx}
\usepackage{dsfont}
\makeatletter
\def\captionof#1#2{{\def\@captype{#1}#2}}
\makeatother
\def\1{\mbox{\bf 1}}
\def\R{\mathbb{R}}

\def\N{\mathbb{N}}
\def\P{\mathbb{P}}
\def\E{\mathbb{E}}
\def\L{\mathbb{L}}

\def\R{\mathbb{R}}

\def\Z{\mathbb{Z}}

\def\c{\mbox{Cov}}

\newtheorem{theo}{Theorem}

\newtheorem{prop}{Proposition}
\newtheorem{Def}{Definition}
\newtheorem{cor}{Corollary}
\newtheorem{Def/Prop}{Definition-Proposition}

\newcounter{exos}
\renewcommand\theexos{\arabic{exos}}

\newcounter{prob}
\renewcommand\theprob{\arabic{prob}}

\begin{document}
\author{Lionel Truquet \footnote{UMR 9194 CREST, ENSAI, Campus de Ker-Lann, rue Blaise Pascal, BP 37203, 35172 Bruz cedex, France. {\it Email: lionel.truquet@ensai.fr}.}
 }
\title{A perturbation analysis of some Markov chains models with time-varying parameters}
\date{}
\maketitle

\begin{abstract}
We study some regularity properties in locally stationary Markov models which 
are fundamental for controlling the bias of nonparametric kernel estimators.
In particular, we provide an alternative to the standard notion of derivative process developed in the literature and that can be used for studying a wide class of Markov processes.
 To this end, for some families of $V-$geometrically ergodic Markov kernels indexed by a real parameter $u$, we give conditions under which the invariant probability distribution is differentiable with respect to $u$, in the sense of signed measures. 
Our results also complete the existing literature for the perturbation analysis of Markov chains, in particular when exponential moments are not finite.
Our conditions are checked on several original examples of locally stationary processes such as integer-valued autoregressive processes, categorical time series or threshold autoregressive processes. 
\end{abstract}

\section{Introduction}
The notion of local stationarity has been introduced in \citet{Dahlhaus} and offers an interesting approach for the modeling of nonstationary time series for which the parameters are continuously changing with the time. In the literature, several stationary models have been extended to a locally stationary version, in particular Markov models defined by autoregressive processes. See for instance \citet{Subba} \citet{Moulines} and \citet{Wu} for linear autoregressive processes, \citet{Dahlhaus2}, \citet{Fryz} and \citet{Truquet2016}
for ARCH processes and a recent contribution of \citet{Dahlhaus3} for nonlinear autoregressive processes. In \citet{Truquet2017}, we have introduced a new notion of local stationarity for general Markov chains models, including most of the autoregressive processes introduced in the references given above but also finite-state Markov chains or integer-valued time series. To define these models, we used time-varying Markov kernels. Let $\left\{Q_u:u\in [0,1]\right\}$ be a family of Markov kernels on the same topological space $\left(E,\mathcal{E}\right)$. We assume that for each $u\in [0,1]$, $Q_u$ has a unique invariant probability measure denoted by $\pi_u$. 
 For an integer $n\geq 1$, we consider $n$ random variables $X_{n,1},X_{n,2},\ldots,X_{n,n}$ such that
\begin{equation}\label{maineqf} 
\P\left(X_{n,t}\in A\vert X_{n,t-1}=x\right)=Q_{t/n}(x,A),\quad (x,A)\in G\times\mathcal{B}(G),\quad 1\leq t\leq n,
\end{equation}
with the convention $X_{n,0}\sim \pi_0$. 
Let us observe that $\left(X_{n,t}\right)_{1\leq t\leq n}$ is a time-inhomogeneous Markov chain as for the locally stationary autoregressive processes of order $1$ introduced in the aforementioned references. 
Then formulation (\ref{maineqf}) is quite general for a locally stationary processes having Markov properties (application to $p-$order Markov process will be also discussed in Section $5$, but as in the homogeneous case, vectorization can be used to get a Markov chain of order $1$). 
The main particularity of our approach, which is similar to that used in the literature of locally stationary processes, is the rescaling by the sample size $n$, taking $Q_{t/n}$ instead of $Q_t$ for the transition kernel at time $t$. The aim of this non standard formulation is to overcome a main drawback of the standard large sample theory, from which it is mainly feasible to estimate parametric models, leading to very arbitrary statistical models for the time-varying Markov kernels $Q_t$. On the other hand, this rescaling allows to use 
a so-called infill asymptotic, from which local inference of some functional parameters defined on the compact unit interval $[0,1]$ remains possible. We defer the reader to the monograph of \citet{Dahl} for a thorough discussion of these asymptotic problems. One of the main issue for making this approach working 
is to show that the triangular array can be approximated marginally (in a sense to precise) by a stationary process with transition kernel $Q_u$ when the ration $t/n$ is close to a point $u\in [0,1]$.   
 
In \citet{Truquet2017}, we proposed a new approach for defining locally stationary Markov chains, using Markov chains techniques. Let us introduce some notations. For two positive integers $t,j$ such that 
$1\leq t\leq n+1-j$, let $\pi^{(n)}_{t,j}$ be the probability distribution of the vector $(X_{n,t},\ldots,X_{n,t+j-1})$ and $\pi_{u,j}$ the probability distribution of the vector $\left(X_1(u),\ldots,X_j(u)\right)$, where $\left(X_t(u)\right)_{t\in\Z}$ denotes a stationary Markov chain with transition kernel $Q_u$. 
Note that $\pi_{u,0}=\pi_u$. In \citet{Truquet2017}, we studied the approximation of $\pi_{t,j}^{(n)}$ by $\pi_{u,j}$ using various probability metrics. One of main idea of the paper is to use contraction/regularity properties for the Markov kernels $Q_u$ which guarantee at the same time such approximation and the decay of some specific mixing coefficients.
We will recall in Section $4$, our approximation result for total variation type norms, from which a large class of locally stationary models can be studied. See also Section $4$ in \citet{Truquet2017} for 
examples of such models and for results on their statistical inference. 

One of the important issues in the statistical inference of locally stationary processes is the curve estimation of some parameters of the kernels $\left\{Q_u: u\in [0,1]\right\}$. However, some parameters of the joint distributions and their regularity, e.g. $\int  fd\pi_u$ for some measurable functionals $f:E\rightarrow \R$, have their own interest for two reasons. 
\begin{enumerate}
\item
First, one can be interested in estimating specific local parameters such as the trend of a time series (which is here the mean of the invariant probability measure) or the local covariance function $u\mapsto \c\left(X_0(u),X_1(u)\right)$. Nonparametric estimation of such functionals typically require to know their regularity, for instance the number of derivatives. For example, estimating the expectation $\int f d\pi_u=\E f\left(X_0(u)\right)$ by a the local linear fit with a kernel density requires the existence of two derivatives for the function $u\mapsto \int f d\pi_u$. See for instance \citet{Fan} for an introduction to local polynomial modeling. We will discuss such a problem in Section $4.3$.
 
\item
Moreover, as discussed in \citet{Truquet2017}, Section $4.5$, when $Q_u(x,dy)=Q_{\theta(u)}(x,dy)$ for a smooth function $\theta:[0,1]\rightarrow\R^d$, getting a bias expression for the local likelihood estimator of $\theta$ requires existence of derivatives for an application of type $u\mapsto \int f d\pi_{2,u}$ where $f:E^2\rightarrow \R$ is a measurable function.
\end{enumerate}
The results stated in \citet{Truquet2017} only guarantee Lipschitz continuity of the applications $u\mapsto \int f d\pi_{u,j}$ for measurable functions $f:E^j\rightarrow \R$. See in particular Proposition $2$ of that paper. One of the aim of the present paper is to complete such results by studying higher-order regularity of such finite-dimensional distributions.

In the recent work of \citet{Dahlhaus3}, the authors study some autoregressive Markov processes with time-varying parameters and defined by iterations of random maps. 
These processes are defined by 
$$X_{n,t}=F_{t/n}\left(X_{n,t-1},\ldots,X_{n,t-p},\varepsilon_t\right),\quad 1\leq t\leq n.$$
Using contraction properties of the random maps $x\mapsto F_u(x,\varepsilon_1)$ in $\L^q-$norms, they study 
the local approximations of $X_{n,t}$ by a stationary process $\left(X_t(u)\right)_{t\in \Z}$ where 
$$X_t(u)=F_u\left(X_{t-1}(u),\ldots,X_{t-p}(u),\varepsilon_t\right),\quad t\in \Z.$$
Differentiability of some functionals of type $u\mapsto \E f\left(X_1(u),\ldots,X_j(u)\right)$ for differentiable functions $f$ are then studied through the notion of a derivative process $dX_t(u)/du$ which is an almost sure derivative of the application $u\mapsto X_t(u)$.  See Proposition $3.8$, Proposition $2.5$ and Theorem $4.8$ in \citet{Dahlhaus3}. The notion of derivative process is fundamental 

Note that here, the process $\left((X_t(u),\ldots,X_{t-p+1}(u))\right)_{t\in\Z}$ form a Markov chain with transition kernel $Q_{u,p}$ defined for $(x_1,\ldots,x_p)\in E^p$ and $(A_1,\ldots,A_p)\in \mathcal{E}^p$ by 
$$Q_{u,p}\left((x_1,\ldots,x_p),A_1\times\cdots\times A_p\right)=\P\left(F_u(x,\varepsilon_0)\in A_p\right)\prod_{i=2}^p\delta_{x_i}(A_{i-1}).$$
 The previous functionals are then defined by some integrals of the invariant probability measure or more generally some integrals of other finite-dimensional distributions of the chain.  Note also that any finite-dimensional distribution of a Markov chain still corresponds to the invariant probability measure of another Markov chain obtained from a vectorization of the initial stochastic process. Studying differentiability properties of an invariant probability measure depending on a parameter is then an important problem.  

For the locally stationary models introduced in \citet{Truquet2017}, the state space is not necessarily continuous, the model is not always defined via contracting random maps and the notion of derivative process is not relevant to evaluate such a regularity. This is in particular the case for count or categorical time series. 
In this paper, our aim is to study directly existence of derivatives for the applications $u\mapsto \pi_{u,j}$ under suitable regularity assumptions for $u\mapsto Q_u$. These derivatives will be understood in the sense of signed measures and using topologies defined by $V-$norms, where $V$ denotes a drift function. See below for further details.
The approach we consider in this paper has two benefits. First, it does not depend on the state space of the Markov process of interest and can be used for lots of locally stationary Markov processes introduced in \citet{Truquet2017} and that cannot be studied using the approach of \citet{Dahlhaus3} (e.g. categorical or count time series). 
Moreover, our approach also applies to the autoregressive processes studied in \citet{Dahlhaus3}. However,  we use Markov chains techniques with small set conditions and stronger regularity assumptions have to be made on the noise distribution.  We defer the reader to the Notes after Proposition \ref{autoreg} for a discussion of the differences between our results and that of \citet{Dahlhaus3} for a time-varying AR$(1)$ process. But, as explained in the same discussion, our results afford a complement to the existing literature because they guarantee differentiability of some maps $u\mapsto \int f d\pi_u$ for non smooth functions $f$ (e.g. the indicator of any Borel set) and allow to consider additional locally stationary autoregressive processes with discontinuous regression functions in space. We also stress that we study differentiability properties of any order whereas \citet{Dahlhaus3} only considered differentiability of order $1$.
The results given in this paper (see in particular Proposition \ref{sufficient2} and Corollary \ref{pMarkov})
are then an alternative to the existing notion of derivative process.

The approach used in this paper has an important connection with the literature of perturbation theory for Markov chains. A central problem in this field is to control an approximation of the invariant probability measure when the Markov kernel of the chain is perturbed. See for instance the recent contribution of \citet{perturb}, motivated by an application to stochastic algorithms. Many works in this field also provide some conditions under which the invariant probability has one or more derivatives with respect to an indexing parameter. See for instance \citet{finite}, \citet{K}, \citet{pflug1}, \citet{vaz} or \citet{Glynn}. For general state spaces, these contributions only focuss on the existence of the first derivative.
Higher-order differentiability is studied using operator techniques in \citet{Heid} or \citet{Heid2}. 
However, as we explain below, these results are restrictive for application to standard time series models.  
Let us first introduce some notations. For a measurable function $V:E\rightarrow [1,\infty)$, we denote by $\mathcal{M}_V(E)$ the set of signed measures $\mu$ on $\left(E,\mathcal{E}\right)$ such that 
$$\Vert\mu\Vert_{V}:=\int V d\vert\mu\vert=\sup_{\vert f\vert\leq V}\int f d\mu<\infty,$$
where $\vert\mu\vert$ denotes the absolute value of the signed measure $\mu$. 
We recall that $\left(\mathcal{M}_V(E),\Vert\cdot\Vert_V\right)$ is a Banach space.
In this paper, we will study differentiability of $u\mapsto \pi_u$, as an application from $[0,1]$ to $\mathcal{M}_V(E)$. The function $V$ will be mainly a drift function for the Markov chain, as in the related references mentioned above. 
We will consider the Markov kernel $Q_u$ as an operator $T_u$ acting on $\mathcal{M}_V(E)$, i.e. $T_u\mu=\mu Q_u$ is the measure defined by
$$\mu Q_u(A)=\int\mu(dx)Q_u(x,A),\quad A\in \mathcal{E}.$$
For a measurable function $g:E\rightarrow \R$ such that $\vert g\vert_V=\sup_{x\in E}\frac{\vert g(x)\vert}{V(x)}<\infty$ , we set $Q_u g(x)=\int Q_u(x,dy)g(y)$. The operator norm of the difference $T_u-T_v$ can be defined by the two following equivalent expressions
$$\Vert T_u-T_v\Vert_{V,V}:=\sup_{\mu\in\mathcal{M}_V(E):\Vert\mu\Vert_V\leq 1}\Vert \mu (P_u-P_v)\Vert_V=\sup_{\vert f\vert_V\leq 1}\left\vert P_u f-P_vf\right\vert_V.$$
Differentiability of the application $u\mapsto \pi_u$, considered as an application form $[0,1]$ to $\mathcal{M}_V(E)$ could be obtained using the results of \citet{Heid} but it is necessary to assume continuity of the application $u\mapsto T_u$ for the previous operator norm. Such continuity assumption is also used in \citet{K}. 
In the literature of perturbation theory, exponential drift functions $V$ are often used and such continuity property can be checked in many examples, such as for some queuing systems considered in \citet{Heid2}.
However, exponential drift functions require exponential moments for the corresponding Markov chain. In time series analysis, existence of exponential moments is a serious restriction.
On the other hand, for power drift functions (another classical choice in the literature of Markov chain), 
this continuity property often fails. For instance, let us consider 
the process $X_t(u)=uX_{t-1}(u)+\varepsilon_t$, $u\in (0,1)$, where $\left(\varepsilon_t\right)_{t\in\Z}$ is a sequence of i.i.d integrable random variables having an absolutely continuous distribution with density $f_{\varepsilon}$. \citet{Hervé} have shown that the corresponding Markov kernel $Q_u(x,dy)=f_{\varepsilon}(y-u x)dy$ is not continuous with respect to $u$, when the classical drift function $V(x)=1+\vert x\vert$ is considered. Additional problems also occur in this example for the derivative operators, obtained by taking the successive derivatives of the conditional density, i.e. $Q_u^{(\ell)}=(-1)^{(\ell)}x^{\ell}f_{\varepsilon}(y-ux)dy$, $\ell=1,2,\ldots$, which are not bounded operators for the operator norm $\Vert\cdot\Vert_{V,V}$.
Boundedness of the derivative operators are required in \citet{Heid} or in \citet{Heid2} for studying the derivatives of $u\mapsto \pi_u$, as an application from $[0,1]$ to $\mathcal{M}_V(E)$. Hence the results 
of the two previous references cannot be applied here. 
For studying differentiability of the invariant probability measure, an alternative result can be found in \citet{HP} (see Appendix A of that paper). This result is applied in \citet{Hervé} to the AR$(1)$ process. However, it is formulated in a very abstract form, using operator theory and its application to on a general class of Markov chain models has not been discussed.  
 
In this paper, we will prove an independent result for studying derivatives of the applications $u\mapsto \pi_u$ or more generally $u\mapsto\pi_{u,j}$ for $j\geq 1$ and that can be applied for a wide class of Markov chains.
This result has some similarities with that of \citet{HP} but our assumptions can be more easily checked and slightly better results can be obtained in the examples we will consider in Section $6$. We defer the reader to the Notes (3.) after Theorem \ref{general} and to the Notes (3.) after Proposition \ref{autoreg} for a discussion. Additionally, for a Markov chain and more generally a $p-$order Markov chain, we provide (see Proposition \ref{sufficient} and Corollary \ref{pMarkov}) easily verifiable conditions on the density of the transition kernels that guarantee differentiability properties for any finite-dimensional distribution of the process.
To our knowledge, the existing literature on the perturbation theory of Markov chains does not contain such 
conditions in a this general context.  
Our approach is particularly useful for models for which some power functions satisfy a drift condition. See Section $4.3$ and Section $5$ for details. 
Moreover, though our results are stated for locally stationary Markov chains, 
one can get a straightforward extension to some parametric models of ergodic Markov processes, using partial derivatives in the multidimensional case. Such modifications will not change the core of our arguments and do not present additional difficulties, we then restrict our study to the case of a parameter $u\in [0,1]$.  

The paper is organized as follows. In Section $2$, we give a general result, formulated using a pure operator-theoretic approach, for getting differentiability properties of an invariant probability measure depending on a parameter. In Section $3$, we give some sufficient conditions on the transition densities of the Markov kernels for applying our result. We also study differentiability of other finite-dimensional distributions of the Markov chain. Section $4$ is devoted to the notion of local stationarity and the control of the bias in kernel smoothing. We also give simple sufficient conditions that ensure both local stationarity and differentiability properties. An extension of our results to $p-$order Markov processes is proposed in Section $5$. We check our assumptions on several examples of locally stationary processes in Section $6$. 
Some of these examples are new or are $p-$order extensions of existing Markov chain models. 
Finally Section $7$ is an Appendix which contains two auxiliary results.

\section{Regularity of an invariant probability with respect to an indexing parameter}

In this section, we consider a family $\left\{P_u: u\in [0,1]\right\}$ of Markov kernels on a topological space $G$ endowed with its Borel $\mathcal{G}=\sigma-$field $\mathcal{B}(G)$. For the locally stationary Markov chains considered in the introduction, we will set $G=E^j$ for $j\geq 1$ and $P_u$ the transition kernel of the Markov chain $Z_t(u)=\left(X_t(u),\ldots,X_{t+j-1}(u)\right)$. See Section $3.3$ for details.
For an integer $k\geq 1$, let $V_0,V_1,\ldots,V_k$ be $k+1$ measurable functions defined on $G$, taking values in $[1,+\infty)$ and such that $V_0\leq V_1\leq \cdots\leq V_k$. For simplicity of notations we set $F_s=\mathcal{M}_{V_s}(G)$ and $\Vert\cdot\Vert_s=\Vert\cdot\Vert_{V_s}$ for $0\leq s\leq k$. We remind that $\left\{\left(F_{\ell},\Vert\cdot\Vert_{\ell}\right): 0\leq \ell\leq k\right\}$ is a family of Banach spaces.
Moreover, $0\leq \ell\leq k-1$, we have $F_{\ell+1}\subset F_{\ell}$ and the injection $$i_{\ell}:\left(F_{\ell+1},\Vert\cdot\Vert_{\ell+1}\right)\rightarrow \left(F_\ell,\Vert\cdot\Vert_\ell\right)$$
is continuous.
For $j=0,1,\ldots,k$, we also denote by $F_{0,j}$ the set of measures $\mu\in F_j$ such that 
$\mu(G)=0$. 
For $0\leq i\leq j\leq k$ and a linear operator $T:\left(F_j,\Vert\cdot\Vert_j\right)\rightarrow \left(F_i,\Vert\cdot\Vert_i\right)$,
we set $\Vert T\Vert_{j,i}=\sup_{\Vert \mu\Vert_j\leq 1}\Vert T \mu\Vert_i$ and $\Vert T\Vert_{0,j,i}=\sup_{\Vert \mu\Vert_j\leq 1, \mu\in F_{0,j}}\Vert T\mu\Vert_i$.
Finally, for each $u\in [0,1]$, we denote by $T_u$ the linear operator acting on the space $F_0$ defined by $T_u \mu=\mu P_u$. For a positive integer $m$, $T_u^m$ will denote the iteration of order $m$ of the operator $T_u$.

\begin{description}
\item
{\bf A1}  We have $T_u F_{\ell}\subset F_{\ell}$ for all $0\leq \ell\leq k$. Moreover, for each $\ell=0,1,\ldots,k$, there exists an integer $m_{\ell}\geq 1$ and a real number $\kappa_{\ell}\in (0,1)$ such that, 
$$\sup_{u\in [0,1]}\Vert T_u^{m_{\ell}}\Vert_{\ell,\ell}\leq \kappa_{\ell},\quad \sup_{u\in [0,1]}\Vert T_u \Vert_{\ell,\ell}<\infty$$
and for each $\mu\in F_{\ell}$, the application $u\mapsto T_u \mu$ is continuous from $[0,1]$ to $\left(F_{\ell},\Vert\cdot\Vert_{\ell}\right)$.

\item
{\bf A2} 
For any $1\leq \ell\leq k$, there exists a continuous linear operator 
$T_u^{(\ell)}:\left(F_{\ell},\Vert\cdot\Vert_{\ell}\right)\rightarrow\left(F_0,\Vert\cdot\Vert_0\right)$ such that for $0\leq s\leq s+\ell\leq k$, $T_u^{(\ell)}F_{s+\ell}\subset F_s$, $\sup_{u\in [0,1]}\left\Vert T_u^{(\ell)}\right\Vert_{s+\ell,s}<\infty$ and for $\mu\in F_{s+\ell}$, the function $u\mapsto T_u^{(\ell-1)}\mu$ 
is differentiable as a function from $[0,1]$ to $F_s$ with continuous derivative $u\mapsto T_u^{(\ell)}\mu$.
We use the convention  $T_u^{(0)}=T_u$.
\item

\end{description}

\begin{theo}\label{general}
Assume the assumptions ${\bf A1-A2}$ hold true. Then the following statements are true.
\begin{enumerate}
\item
\begin{itemize}
\item
For each $u\in [0,1]$, the operator $I-T_u$ defines an isomorphism on each space $\left(F_{0,\ell}\,\Vert\cdot\Vert_{\ell}\right)$ for $0\leq \ell\leq k$. Moreover the inverse of $I-T_u$ is given by $(I-T_u)^{-1}=\sum_{k\geq 0}T_u^k$. 
\item
We have $\max_{0\leq\ell\leq k}\sup_{u\in[0,1]}\left\Vert (I-T_u)^{-1}\right\Vert_{\ell,\ell}<\infty$. 
\item
For $0\leq \ell\leq k$ and $\mu\in F_{0,\ell}$, the application $u\mapsto (I-T_u)^{-1}\mu$ is continuous as an application from $[0,1]$ to $F_{\ell}$. 
\item
Moreover, for each $u\in [0,1]$, we have for $0\leq \ell\leq k-1$ and $\mu\in F_{\ell+1}$,
$$\lim_{h\rightarrow 0}\left\Vert \frac{(I-T_{u+h})^{-1}\mu-(I-T_u)^{-1}\mu}{h}-(I-T_u)^{-1}T_u^{(1)}(I-T_u)^{-1}\mu\right\Vert_{\ell}=0.$$
\end{itemize}
\item
For each $u\in [0,1]$, there exists a unique probability measure $\mu_u$ such that $T_u\mu_u=\mu_u$ ($\mu_u$ is an invariant probability for $P_u$). Moreover $\mu_u\in F_k$.
\item
The application $f:[0,1]\rightarrow F_k$ defined by $f(u)=\mu_u$, for $u\in [0,1]$, is continuous.
Moreover there exist some functions $f^{(0)},\ldots,f^{(k)}$ such that $f^{(0)}=f$ and 
\begin{itemize}
\item
for $1\leq \ell\leq k$, the application $f^{(\ell)}:[0,1]\rightarrow F_{0,k-\ell}$ is continuous,
\item
for $1\leq\ell\leq k$ and $u\in[0,1]$,  $\lim_{h\rightarrow 0}\left\Vert \frac{f^{(\ell-1)}(u+h)-f^{(\ell-1)}(u)}{h}-f^{(\ell)}(u)\right\Vert_{k-\ell}=0,$ 
\item
the derivatives of $f$ are given recursively by 
$$f^{(\ell)}(u)=\sum_{s=1}^{\ell}\begin{pmatrix}\ell\\s\end{pmatrix}\left(I-T_u\right)^{-1} T_u^{(s)}f^{(\ell-s)}(u).$$
\end{itemize}
\end{enumerate}
\end{theo}

\paragraph{Notes} 
\begin{enumerate}
\item
When $V_0=V_1=\cdots=V_k=V$, existence of the derivatives for the invariant probability measures is studied in \citet{Heid}. 
One can show that the condition {\bf $C^k$} used for stating their result entails {\bf A2} because they use a continuity assumption of the derivative operators with respect to the $V-$operator norm. On the other hand, their geometric ergodicity result (see Result $2$ in their paper) for each kernel $P_u$, the measure and the continuity assumption of the kernel for the $V-$operator norm entails the contraction {\bf A1} (for the contraction coefficient, see Section $3.2$ below).
We also deduce from our result the following Taylor-Lagrange formula that will be useful for controlling the bias of kernel estimators in Section $4.2$. For $u\in [0,1]$ and $h\in\R$ such that $u+h\in [0,1]$, set 
$M=\sup_{v\in [0,1]}\Vert f^{(k)}(v)\Vert_0$. We then have
\begin{equation}\label{super}
\Vert f(u+h)-f(u)-\sum_{\ell=1}^{k-1}\frac{f^{(\ell)}(u)}{\ell!}h^{\ell}\Vert_0\leq \frac{M\vert h\vert^k}{k!}.
\end{equation}

\item
Let us discuss our assumptions. Assumption {\bf A1} guarantees the stability of the spaces $\mathcal{M}_{V_s}(G)$ by the application $T_u$ (i.e. $\mu\in \mathcal{M}_{V_s}(G)\rightarrow \mu P_u\in \mathcal{M}_{V_s}(G)$). 
The contraction condition in the second part of this assumption guarantees some invertibility properties of the operator $I-T_u$ (see point $1$ of Theorem \ref{general}) that are needed for getting an expression of the derivatives of $u\mapsto \mu_u$. 
One can see that our assumptions involve some measure spaces of increasing regularity $\mathcal{M}_{V_k}(G)\subset \cdots\subset \mathcal{M}_{V_0}(G)$. 
Assumption {\bf A2} allows the derivative operators of the Markov kernel to be only bounded for an operator norm involving a weaker final topology. This is particularly useful when the derivatives operators do not preserve a measure space of given regularity. For instance, for the AR(1) process $X_k(u)=a(u)X_{k-1}(u)+\varepsilon_k$ with a noise density $f_{\varepsilon}$, we have $T_u\mu(dy)=\int \mu(dx)f_{\varepsilon}(y-a(u)x)dy$ and a natural candidate for $T_u^{(\ell)}$ is 
$$T_u^{(\ell)}\mu(dy)=a^{(\ell)}(u)\int \mu(dx)(-x)^{\ell}f^{(\ell)}_{\varepsilon}(y-a(u)x)dy.$$ 
Setting $V_s=1+\vert x\vert^s$ , one can see that $\left\vert T_u^{(\ell)}\mu\right\vert\cdot V_s\leq C \vert\mu\vert\cdot V_{s+\ell}$ for a positive constant $C$. This means that $\mu$ has to have a moment of order $s+\ell$ for getting a finite upper bound in the previous inequality.
This problem does not occur on this example when the $V_s'$s are some exponential functions and the noise density and its derivatives have exponential moments. See in particular Proposition \ref{expo} given in the Appendix. However, we do not want to use this restrictive moment condition.  

\item
The idea of introducing spaces of increasing regularity (as $\mathcal{M}_{V_k}(G)\subset\cdots\subset \mathcal{M}_{V_0}(G)$ in our result) can also be found in \citet{HP} (see Annex $A$ of that paper). In Proposition $A$ of that paper, the authors study regularity properties of some resolvent operators depending on a parameter. They also used an operator theoretic approach. An application of this result to study the regularity of the invariant probability measure of an AR$(1)$ process with respect to its autoregressive coefficient is given in \citet{Hervé}, Proposition $1$. However, application of such a result requires in our context to introduce additional operator norms for getting continuity properties of applications $u\mapsto T_u^{(\ell)}$, as applications form $[0,1]$ to some spaces of linear operators. See in particular the proof of Proposition $1$ in \citet{Hervé} and the transformation $T_0$ introduced in the proof of their Lemma $1$. Here, in {\bf A2}, we prefer to use pointwise continuity/differentiability assumptions for some applications $u\mapsto T_u^{(\ell)}\mu$ and that are sufficient for getting our result. We found our formulation easier to understand.
We also defer the reader to the Notes (3.) after Proposition \ref{autoreg} for a comparison of our result with that of \citet{Hervé} for an AR$(1)$ process. 

\item
In Assumption {\bf A2}, we assume that the operators $T_u^{(\ell)}$ satisfies some kind of weak continuity 
or weak differentiability with respect to $u$, in the sense that continuity and differentiability do not not hold for operator norms but simply for some applications $u\mapsto T_u^{(\ell)}\mu$. In the literature of perturbation of Markov chains, a notion of weak continuity or differentiability for measures depending on parameters can be found in \citet{pflug2} (see Section $3.2$). Our condition is stronger since for an individual measure $\mu$, the application $u\mapsto \mu P_u f$ is required to be continuous or differentiable but uniformly over a class of functions $f$. In contrast, \citet{pflug2} defined these notions for a fixed function $f$. But note that our final result entails existence of derivatives for the topology defined by some $V-$norms, which is stronger than getting derivatives for $u\mapsto \int fd\mu_u$ for a single function $f$.       
\end{enumerate}

\paragraph{Proof of Theorem \ref{general}}

\begin{enumerate}
\item
\begin{itemize}
\item
First, one can note that $\left(F_{0,\ell},\Vert\cdot\Vert_{\ell}\right)$ is a closed vector subspace of $\left(F_{\ell},\Vert\cdot\Vert_{\ell}\right)$ and then a Banach space. Moreover, From Assumption {\bf A1}, the series $\sum_{k\geq 0} T_u^k$, considered as an operator from $F_{0,\ell}$ to $F_{0,\ell}$ is normally convergent for the norm $\Vert\cdot\Vert_{0,\ell,\ell}$ and is the inverse of $I-T_u$.
Then $I-T_u$ defines an isomorphism on the space $\left(F_{0,\ell},\Vert\cdot\Vert_{\ell}\right)$. 
\item
Using the expression $(I-T_u)^{-1}=\sum_{k\geq 0}T_u^k$, the second assertion is a consequence of Assumption {\bf A1}.
\item
Next, we show that for $0\leq\ell\leq k$ and $\mu\in F_{0,\ell}$, the application $u\mapsto (I-T_u)^{-1} \mu$ is continuous as an application from $[0,1]$ to $F_{0,\ell}$. Considering all the operators as operators from $F_{0,\ell}$ to $F_{0,\ell}$, we use the decomposition
\begin{equation}\label{decdecdec} 
(I-T_{u+h})^{-1}-(I-T_u)^{-1}=(I-T_{u+h})^{-1}(T_{u+h}-T_u)(I-T_u)^{-1}.
\end{equation}
From the previous point, we have $\sup_{u\in [0,1]}\Vert (I-T_u)^{-1}\Vert_{0,\ell,\ell}<\infty$ and $(I-T_u)^{-1}\mu$ is an element of $F_{0,\ell}$. Moreover, if $\nu\in F_{\ell}$, Assumption {\bf A1} guarantees the continuity of the application $v\mapsto T_v\nu$ as an application from $[0,1]$ to $F_{\ell}$.
Using (\ref{decdecdec}), the continuity of the application $u\mapsto (I-T_u)^{-1} \mu$ follows.
\item
Finally, if $\mu\in F_{0,\ell+1}$, we show that the application $u\mapsto (I-T_u)^{-1}\mu$ is differentiable
as an application from $[0,1]$ to $F_{0,\ell}$.
Setting $z_{u,h}=h^{-1}\left(T_{u+h}-T_u\right)(I-T_u)^{-1}\mu$, we deduce from Assumption {\bf A2} that $\lim_{h\rightarrow 0}z_{u,h}=z_u=T_u^{(1)}(I-T_u)^{-1}\mu$ in $\left(F_{0,\ell},\Vert\cdot\Vert_{\ell}\right)$. 
We use the decomposition 
\begin{eqnarray*}
h^{-1}\left[(I-T_{u+h})^{-1}x-(I-T_u)^{-1}x\right]&=&(I-T_{u+h})^{-1}z_{u,h}\\
&=& (I-T_{u+h})^{-1}(z_{u,h}-z_u)+(I-T_{u+h})^{-1}z_u.
\end{eqnarray*}
From the previous point, we have $\lim_{h\rightarrow 0}(I-T_{u+h})^{-1}z_u=(I-T_u)^{-1}z_u$ in $\left(F_{0,\ell},\Vert\cdot\Vert_{\ell}\right)$.
Moreover,
$$\left\Vert (I-T_{u+h})^{-1}(z_{u,h}-z_u)\right\Vert_{\ell}\leq \sup_{u\in [0,1]}\left\Vert (I-T_u)^{-1}\right\Vert_{\ell,\ell}\Vert z_{u,h}-z_u\Vert_{\ell}\stackrel{h\rightarrow 0}{\rightarrow}0.$$
This shows that the application $u\mapsto (I-T_u)^{-1}\mu$ is differentiable, as an application from $[0,1]$ to $F_{0,\ell}$, with derivative $u\mapsto (I-T_u)^{-1}T_u^{(1)}(I-T_u)^{-1}\mu$.
\end{itemize} 
\item
The space $F_{k,1}=\left\{\mu\in F_k: \mu \mbox{ is a probability measure }\right\}$ endowed with the norm $\Vert\cdot\Vert_k$ is a complete metric space. From Assumption {\bf A1} and the fixed point theorem, there exists a unique probability measure $\mu_u$ in $F_k$ such that $\mu_u P_u=\mu_u$. But $\mu_u$ is in fact the single invariant probability measure for $P_u$. Indeed, since for any $x\in G$, we have $\delta_x\in F_{k,1}$, 
the fixed point theorem applied in $F_{k,1}$ entails that $\lim_{n\rightarrow \infty}\int f(y)P_u^{n}(x,dy)=\int f(y)\mu_u(dy)$ for all measurable function $f:G\rightarrow \R$ bounded by one. 
If $\overline{\mu}_u$ is an invariant probability measure, we get
from the Lebesgue theorem, 
$$\int f(y)\overline{\mu}_u(dy)=\int f(y)P_u^n(x,dy)\overline{\mu}(dx)\rightarrow_{n\rightarrow \infty}\int f(y)\mu_u(dy).$$
Necessarily, $\overline{\mu}_u=\mu_u$ and $\mu_u$ is then the unique invariant probability measure for $P_u$.    
\item
We first show that $f$ is continuous. 
We have $f(u+h)-f(u)=(I-T_{u+h})^{-1}(T_{u+h}-T_u)f(u)$. From Assumption {\bf A1}, we have 
$\lim_{h\rightarrow 0}\left\Vert T_{u+h}f(u)-T_u f(u)\right\Vert_k=0$. Note that $(T_{u+h}-T_u)f(u)$ is an element of $F_{0,k}$. Using the second assertion of point $1.$ of the theorem, we get $\lim_{h\rightarrow 0}\left(f(u+h)-f(u)\right)=0$ in $\left(F_{0,k},\Vert\cdot\Vert_k\right)$.
 
Next, we prove the existence of the derivatives and their properties by induction on $\ell$ with $1\leq \ell\leq k$. 
\begin{enumerate}
\item
First, we assume that $\ell=1$.  Using the same decomposition as for proving continuity of $f$, we have 
\begin{eqnarray*}
&&\frac{f(u+h)-f(u)}{h}\\
&=&\left(I-T_{u+h}\right)^{-1}\frac{T_{u+h}-T_u}{h}\mu_u\\
&=& \left(I-T_{u+h}\right)^{-1}\left[\frac{T_{u+h}-T_u}{h}\mu_u-T_u^{(1)}\mu_u\right]+\left(I-T_{u+h}\right)^{-1}T_u^{(1)}\mu_u.
\end{eqnarray*}
Here we consider the operators $T_{u+h}-T_u$ and $T_u^{(1)}$ as operator from $F_k$ to $F_{0,k-1}$. The operators $(I-T_u)^{-1}$, $u\in [0,1]$, are considered as operators from $F_{0,k-1}$ to $F_{0,k-1}$. 
From Assumption {\bf A2}, we have 
$$\lim_{h\rightarrow 0}\left\Vert \frac{T_{u+h}\mu_u-T_u \mu_u}{h}-T^{(1)}_u \mu_u\right\Vert_{k-1}=0.$$
From the second and the third assertions of the point $1.$, we get 
$$\lim_{h\rightarrow 0}\left\Vert \frac{f(u+h)-f(u)}{h}-f^{(1)}(u)\right\Vert_{k-1}=0,$$
where $f^{(1)}(u)=(I-T_u)^{-1}T_u^{(1)}\mu_u$. 
It remains to prove the continuity of $f^{(1)}$ as an application from $[0,1]$ to $F_{k-1}$. 
As previously, it is sufficient to show that 
$$\lim_{h\rightarrow 0}\left\Vert T_{u+h}^{(1)}\mu_{u+h}-T_u^{(1)}\mu_u\right\Vert_{k-1}=0.$$
But this is a consequence of the continuity of $f$ and of Assumption {\bf A2}, using the decomposition 
$$T_{u+h}^{(1)}\mu_{u+h}-T_u^{(1)}\mu_u=\left[T_{u+h}^{(1)}-T_u^{(1)}\right]\mu_u+T_{u+h}^{(1)}\left[\mu_{u+h}-\mu_u\right].$$
This shows the result for $\ell=1$.

\item
Now let us assume that for $1\leq\ell\leq k-1$, $f$ has $\ell$ derivatives such that for $1\leq s\leq \ell$ and $u\in[0,1]$, 
the function $f^{(s)}:[0,1]\rightarrow F_{0,k-s}$ is  continuous, 
$$\lim_{h\rightarrow 0}\left\Vert \frac{f^{(s-1)}(u+h)-f^{(s-1)}(u)}{h}-f^{(s)}(u)\right\Vert_{k-s}=0$$ 
and  
$$f^{(\ell)}(u)=\sum_{s=1}^{\ell}\begin{pmatrix}\ell\\s\end{pmatrix}\left(I-T_u\right)^{-1} T_u^{(s)}f^{(\ell-s)}(u).$$
\begin{itemize}
\item
For $1\leq s\leq \ell$, we set $z_u=T_u^{(s)}f^{(\ell-s)}(u)$ and we consider $T_u^{(s)}$ as an operator from $F_{k-\ell+s}$ to $F_{k-\ell}$. We are going to show that the application $u\mapsto z_u$ from $[0,1]$ to $F_{0,k-\ell}$ has a derivative. We have 
$$\frac{z_{u+h}-z_u}{h}=\frac{T^{(s)}_{u+h}-T_u^{(s)}}{h}f^{(\ell-s)}(u)+T_{u+h}^{(s)}\frac{f^{(\ell-s)}(u+h)-f^{(\ell-s)}(u)}{h}.$$
Since $f^{(\ell-s)}(u)\in F_{k-\ell+s}$, we have from Assumption {\bf A2}, 
$$\lim_{h\rightarrow 0}\left\Vert \frac{T^{(s)}_{u+h}-T_u^{(s)}}{h}f^{(\ell-s)}(u)-T^{(s+1)}f^{(\ell-s)}(u)\right\Vert_{k-\ell-1}=0.$$
Next we set $w_{u,h}=\frac{f^{(\ell-s)}(u+h)-f^{(\ell-s)}(u)}{h}$. By the induction hypothesis, we have 
$$\lim_{h\rightarrow 0}\left\Vert w_{u,h}-f^{(\ell-s+1)}(u)\right\Vert_{k-\ell+s-1}=0.$$
Using Assumption {\bf A2}, we have $\sup_{u\in [0,1]}\left\Vert T_u^{(s)}\right\Vert_{k-\ell+s-1,k-\ell-1}<\infty$. 
Then we get 
$$\lim_{h\rightarrow 0}\left\Vert T_{u+h}^{(s)}\left(w_{u,h}-f^{(\ell-s+1)}(u)\right)\right\Vert_{k-\ell-1}=0.$$
Using again Assumption {\bf A2}, we have 
$$\lim_{h\rightarrow 0}\left\Vert T_{u+h}^{(s)}f^{(\ell-s+1)}(u)-T_u^{(s)}f^{(\ell-s+1)}(u)\right\Vert_{k-\ell-1}=0.$$
This shows that 
$$\lim_{h\rightarrow 0}\left\Vert \frac{z_{u+h}-z_u}{h}-T_u^{(s+1)}f^{(\ell-s)}(u)-T_u^{(s)}f^{(\ell-s+1)}(u)\right\Vert_{k-\ell-1}=0.$$
In the sequel we set $z_u^{(1)}=T_u^{(s+1)}f^{(\ell-s)}(u)+T_u^{(s)}f^{(\ell-s+1)}(u)$. 
\item
Next we compute the derivative of $u\mapsto y_u=(I-T_u)^{-1}z_u$, as an application from $[0,1]$ to $F_{0,k-\ell-1}$. 
We have 
$$\frac{y_{u+h}-y_u}{h}=\frac{(I-T_{u+h})^{-1}-(I-T_u)^{-1}}{h}z_u+(I-T_{u+h})^{-1}\left(\frac{z_{u+h}-z_u}{h}-z_u^{(1)}\right)+(I-T_{u+h})^{-1} z_u^{(1)}.$$
Using Assumption {\bf A2} and some previous results, we get 
$$\lim_{h\rightarrow 0}\left\Vert \frac{y_{u+h}-y_u}{h}-(I-T_u)^{-1}T_u^{(1)}(I-T_u)^{-1}z_u-(I-T_u)^{-1}z_u^{(1)}\right\Vert_{k-\ell-1}=0.$$
In the sequel, we set 
$$t^{(\ell,s)}(u)=(I-T_u)^{-1}T_u^{(1)}(I-T_u)^{-1}z_u+(I-T_u)^{-1}z_u^{(1)}.$$
\item
Finally we get in $\left(F_{k-\ell-1},\Vert\cdot\Vert_{k-\ell-1}\right)$,
$$\lim_{h\rightarrow 0}\frac{f^{(\ell)}(u+h)-f^{(\ell)}(u)}{h}=f^{(\ell+1)}(u),$$
where 
$$f^{(\ell+1)}(u)=\sum_{s=1}^{\ell}\begin{pmatrix}\ell\\ s\end{pmatrix}
t_u^{(\ell,s)}.$$
The expression for $f^{(\ell+1)}(u)$ given in the statement of the theorem follows from straightforward computations.
\item
Finally, using the induction hypothesis, the function $f^{(\ell+1-s)}$ is continuous as an application from $[0,1]$ to $F_{k-\ell+s-1}$, for each $1\leq s\leq \ell+1$. The proof of the continuity of $f^{(\ell+1)}$ is then similar to the proof of the continuity of $f^{(1)}$. 
\end{itemize}

The properties of the successive derivatives $f^{(1)},\ldots,f^{(k)}$ follow by induction and the proof of Theorem \ref{general} is now complete.$\square$
\end{enumerate}
 \end{enumerate}

\section{Sufficient conditions}

We now provide some sufficient conditions for ${\bf A1-A2}$. The two previous assumptions are not easy to verify and we want to provide some conditions that can be more easily checked for practical examples. We also provide a natural expression for the derivative operators $T_u^{(\ell)}$.
Assumption {\bf B1} given below is related to uniform ergodicity. Since there are several ways of checking this new assumption, we discuss it in Section $3.2$. 
 
In what follows, we assume that the kernel $P_u$ is defined by
$$P_u(x,A)=\int_A f(u,x,y)\gamma(x,dy),\quad A\in \mathcal{B}(G),$$
where $f:[0,1]\times G^2\rightarrow \R_+$ is a measurable function and $\gamma$ is a kernel not depending on $u$.

\subsection{A sufficient set of conditions}

In order to check ${\bf A1-A2}$, we make some regularity assumptions on the family of conditional densities $\left\{f(u,\cdot,\cdot): u\in [0,1]\right\}$. Let $k$ be a positive integer and $V_k\geq V_{k-1}\geq\cdots \geq V_0$ some measurable applications from $G$ to $[1,\infty)$ such that the following conditions are satisfied. 
\begin{description}
\item[B1]
For $\ell=0,1,\ldots, k$, the family of Markov kernels $\left\{P_u:u\in [0,1]\right\}$ is simultaneously $V_{\ell}-$uniformly ergodic, i.e there exists $\kappa_{\ell}\in (0,1)$ such that,
$$\sup_{u\in [0,1]}\sup_{x\in G}\frac{\Vert \delta_x P_u^n-\mu_u\Vert_{\ell}}{V_{\ell}(x)}=O\left(\kappa_{\ell}^n\right),$$
where the unique invariant probability measure $\mu_u$ of $P_u$ satisfies $\mu_u V_k<\infty$.

\item[B2] For all $(x,y)\in G^2$, the function $u\mapsto f(u,x,y)$ is $k-$times continuously differentiable and for $1\leq \ell\leq k$, we denote by $\partial^{(\ell)}_1f$ its partial derivative of order $\ell$.

\item[B3]
There exist $C>0$ such that for integers $0\leq s\leq s+\ell\leq k$ and $x\in G$, 
\begin{equation}\label{Hyp1}
\sup_{u\in [0,1]}\int \left\vert \partial^{(\ell)}_1f(u,x,y)\right\vert V_s(y)\gamma(x,dy)\leq C V_{s+\ell}(x)
\end{equation}
and for each $u\in [0,1]$, 
\begin{equation}\label{Hyp2}
\lim_{h\rightarrow 0}\int \left\vert \partial^{(k-s)}_1f(u+h,x,y)-\partial^{(k-s)}_1f(u,x,y)\right\vert V_s(y)\gamma(x,dy)=0.
\end{equation}
\end{description}

\begin{cor}\label{sufficient}
The assumptions {\bf B1-B3} entail the assumptions {\bf A1-A2}. Moreover the conclusions of Theorem \ref{general} are valid for the derivative operators 
$$T_u^{(\ell)}\mu=\int \mu(dx)\partial^{(\ell)}_1f(u,x,y)\gamma(x,dy),\quad 1\leq \ell\leq k,\quad \mu\in \mathcal{M}_{V_{\ell}}(G).$$
\end{cor}

\paragraph{Proof of Corollary \ref{sufficient}}
\begin{enumerate}
\item
We first check {\bf A1}.
If $P$ is a Markov kernel on $\left(G,\mathcal{B}(G)\right)$, we define the following Dobrushin contraction coefficient 
\begin{equation}\label{Dobcont}
\Delta_V(P):=\sup_{\mu\in\mathcal{M}_V(G),\mu\neq 0,\mu(G)=0}\frac{\Vert \mu P\Vert_V}{\Vert\mu\Vert_V}=\sup_{x,y\in G,x\neq y}\frac{\Vert \delta_xP-\delta_yP\Vert_V}{V(x)+V(y)}.
\end{equation}
See for instance \citet{MoulDouc}, Lemma $6.18$ for the second expression. 
Note also that, with the notations of Section $2$, we have if $T\mu=\mu P$, $\Vert T\Vert_{0,\ell,\ell}=\Delta_{V_{\ell}}(P)$.

First, note that from (\ref{Hyp1}) applied with $\ell=0$, we have $T_u F_s\subset F_s$ and $\sup_{u\in[0,1]}\Vert T_u\Vert_{s,s}<\infty$ for $s=0,1,\ldots,k$. 
Moreover, we have the bound 
\begin{equation}\label{banddob}
\Vert T_u^n\Vert_{0,\ell,\ell}=\Delta_{V_{\ell}}\left(P_u^n\right)\leq \sup_{x\in G}\frac{\Vert \delta_x P_u^n-\pi_u\Vert_{\ell}}{V_{\ell}(x)}.
\end{equation}
This bound can be found for instance in \citet{perturb}, Lemma $3.2$. 
For completeness, we repeat the argument.
We have, using the inequality $(a+b)/(c+d)\leq \max\{a/c,b/d\}$ valid for all positive real numbers $a,b,c,d$,
$$\Delta_V\left(P_u^n\right)\leq \sup_{x\neq y}\frac{\Vert\delta_x P_u^n-\delta_yP_u^n\Vert_V}{V(x)+V(y)}\leq \sup_{x\neq y}\frac{\Vert \delta_x P_u^n-\mu_u\Vert_V+\Vert\delta_y P_u^n-\mu_u\Vert_V}{V(x)+V(y)}\leq \sup_{x\in G}\frac{\Vert \delta_x P_u^n-\mu_u\Vert_V}{V(x)},$$
which shows (\ref{banddob}).
This entails the existence of an integer $m_{\ell}\geq 1$ such that $\sup_{u\in [0,1]}\Vert T_u^{m_{\ell}}\Vert_{0,\ell,\ell}<1$. 
It remains to show that if $\mu \in F_{\ell}$, $u\mapsto T_u\mu$ is continuous, as an application from $[0,1]$ to $F_{\ell}$. We have 
$$\Vert T_{u+h}\mu-T_u\mu\Vert_{\ell}\leq \int\int \vert\mu\vert(dx)\gamma(x,dy)V_{\ell}(y)\left\vert f(u+h,x,y)-f(u,x,y)\right\vert.$$
We will use the Lebesgue theorem. Using the inequality $V_{\ell}\leq V_k$ and Assumption {\bf B3} (\ref{Hyp2})with $s=k$,
$$\lim_{h\rightarrow 0}c_h(u,x):=\int\gamma(x,dy)V_{\ell}(y)\left\vert f(u+h,x,y)-f(u,x,y)\right\vert=0,\quad x\in G.$$
Moreover, from {\bf B3} (\ref{Hyp1}) applied to the derivative of order $0$, we have
$c_h(u,x)\leq 2CV_{\ell}(x)$ and $V_{\ell}$ is $\vert\mu\vert-$integrable. The Lebesgue theorem then applies and gives $\lim_{h\rightarrow 0}T_{u+h}\mu=T_u\mu$ in $F_{\ell}$ and the last assertion in {\bf A1} follows.  
\item
Next, we check the assumption {\bf A2}. We first notice that
for $0\leq s\leq s+\ell\leq k$ and $\mu\in F_{s+\ell}$, we have from {\bf B3} (\ref{Hyp1}), 
\begin{eqnarray*}
\Vert T_u^{(\ell)}\mu\Vert_s&\leq& \int\vert\mu\vert(dx)\int \gamma(x,dy)\left\vert\partial^{(\ell)}_1f(u,x,y)\right\vert V_s(y)\\
&\leq& C \int\vert\mu\vert(dx)V_{s+\ell}(x)=C\Vert\mu\Vert_{s+\ell}.
\end{eqnarray*}
This shows that $T_u^{(\ell)}F_{s+\ell}\subset F_s$ and $\sup_{u\in[0,1]}\Vert T_u^{(\ell)}\Vert_{s+\ell,s}\leq C$. 
Next, for $\mu\in F_s$, we show the continuity of the application $u\mapsto T_u^{(\ell)}\mu$, as an application from $[0,1]$ to $F_{s+\ell}$.
We have 
$$\Vert T_{u+h}^{(\ell)}\mu-T_u^{(\ell)}\mu\Vert_s\leq \int \vert\mu\vert(dx)\int \mu_0(x,dy)\left\vert\partial_1^{(\ell)}f(u+h,x,y)
-\partial_1^{(\ell)}f(u,x,y)\right\vert V_s(y).$$
From the assertion (\ref{Hyp1}) in {\bf B3} and the Lebesgue theorem, it is enough to prove that for all $x\in G$,
$$\lim_{h\rightarrow 0}\int \gamma(x,dy)\left\vert\partial_1^{(\ell)}f(u+h,x,y)
-\partial_1^{(\ell)}f(u,x,y)\right\vert V_s(y)=0.$$
We consider two cases.
\begin{itemize}
\item
If $s+\ell=k$, this continuity is a direct consequence of the assertion (\ref{Hyp2}) of Assumption {\bf B3}. 
\item
We next assume that $s+\ell+1\leq k$. We have
$$\int \gamma(x,dy)\left\vert\partial_1^{(\ell)}f(u+h,x,y)
-\partial_1^{(\ell)}f(u,x,y)\right\vert V_s(y)
\leq h\sup_{v\in[0,1]}\int\gamma(x,dy) \left\vert\partial_1^{(\ell+1)}f(v,x,y)\right\vert V_s(y).$$
Then the result follows from the assumption {\bf B3} (\ref{Hyp1}).
\end{itemize}

Finally, we show the differentiability property of the operators. For $\mu\in F_{s+\ell}$, we have, using the mean value theorem, 
\begin{eqnarray*}
&&\Vert \frac{T^{(\ell-1)}_{u+h}\mu-T_u^{(\ell-1)}\mu}{h}-T_u^{(\ell)}\mu\Vert_s\\
&\leq& \sup_{v\in [u,u+h]}\int\vert\mu\vert(dx)\int \gamma(x,dy)\left\vert\partial_1^{(\ell)}f(v,x,y)-\partial_1^{(\ell)}f(u,x,y)\right\vert V_s(y).
\end{eqnarray*}
The result follows by using the same arguments as in the proof of the continuity of the application $u\mapsto T^{(\ell)}_u\mu$. 
This completes the proof of Corollary \ref{sufficient}.$\square$ 
\end{enumerate}

\paragraph{Note.} When $\phi:G\rightarrow [1,\infty)$ is a measurable function such that for some $d\leq d_0$, $0\leq \ell\leq k$ and $q_0,q_1,\ldots,q_k>0$,
$$\int \gamma(x,dy)\left\vert \partial_1^{(\ell)}f(u,x,y)\right\vert \phi(y)^d\leq C \phi(x)^{d+q_{\ell}},$$
assumption {\bf B3} (\ref{Hyp1}) is checked by setting $V_{\ell}(x)=\phi(x)^{d+q\ell}$ with $q=\max\left(q_1,q_2/2,\ldots,q_k/k\right)$ and assuming that $d+qk\leq d_0$. 

\subsection{Simultaneous uniform ergodicity}\label{explain}

Assumption {\bf B1} is related to a simultaneous $V-$uniform ergodicity condition. Let us first give a precise definition of this notion.
\begin{Def}\label{manque}
We will say that a family of Markov kernel $\left\{P_u:u\in [0,1]\right\}$ satisfies a simultaneous $V-$uniform ergodicity condition if there exists $C>0$ and $\kappa\in(0,1)$ such that for all $u\in [0,1]$ and all $x\in G$, 
$$\Vert \delta_x P_u^n-\mu_u\Vert_{V}\leq C V(x)\kappa^n.$$
\end{Def}
This notion plays a central rule in our results and it is then important to provide sufficient conditions for ${\bf B1}$. We also point out that this notion of simultaneous uniform ergodicity replaces stronger assumptions made in \citet{Heid}. These authors used pointwise uniform ergodicity and a continuity property for the application $u\mapsto P_u$, in the sense that 
$$\lim_{h\rightarrow 0}\Vert P_{u+h}-P_u\Vert_{V,V}=0.$$
See in particular Definition $3$ and Condition $1-4$ of that paper. For simplicity, we justify why the two previous conditions imply simultaneous uniform ergodicity in a separate result (see Proposition \ref{compar}   in Section \ref{Appendix}).

For a single Markov kernel, $V-$uniform ergodicity is generally obtained under a drift condition and a small set condition. See \citet{MT}, Chapter $16$ for details. Let us first recall the definition of a small set.  For a Markov kernel $P$ on $\left(G,\mathcal{B}(G)\right)$, a set $C\in\mathcal{B}(G)$ is called a $(\eta,\nu)-$small set, where $\eta$ a positive real number and $\nu$ a probability measure on $\left(G,\mathcal{B}(G)\right)$ if  
$$P(x,A)\geq \eta \nu(A),\mbox{  for all } A\in \mathcal{B}(G) \mbox{ and all } x\in C.$$
 
We now present two approaches for getting simultaneous uniform ergodicity.

\subsubsection{Simultaneous $V-$uniform ergodicity via drift and small set conditions}\label{rat}
When simultaneous drift and small set conditions are satisfied, a result of \citet{HM1} can be used to check simultaneous $V-$uniform ergodicity. 
For simplicity we introduce the following condition. For $\lambda\in (0,1)$, $b,\eta, r>0$ and $\nu$ a probability measure on $\left(G,\mathcal{B}(G)\right)$, 
we will say that a Markov kernel $P$ satisfies the condition  
$\mathcal{C}\left(V,\lambda,b,r,\eta,\nu\right)$ if 
\begin{equation}\label{centplus}  
P V\leq \lambda V+b\quad \mbox{and}\quad \left\{x\in G :V(x)\leq r\right\}\mbox{ is a } (\eta,\nu)\mbox{ small set}.
\end{equation}
If there exists an integer $m\geq 1$ such that all the Markov kernels $P_u^m$, $u\in [0,1]$, satisfy the condition $\mathcal{C}\left(V,\lambda,b,r,\eta,\nu\right)$ for $r>\frac{2b}{1-\lambda}$ 
and if there exists $K>0$ such that $P_uV \leq KV$ for all $u\in [0,1]$, Theorem $1.3$ in \citet{HM1} guaranty the existence of $\alpha\in (0,1)$ and $\delta>0$, not depending on $u\in [0,1]$ such that $\Delta_{V_{\delta}}\left(P_u^m\right)\leq \alpha$ with $V_{\delta}=1+\delta V$ (see (\ref{Dobcont}) for the definition of $\Delta_V$). Actually, the result of \citet{HM1} is stated for a single Markov kernel but inspection of the proof shows that the coefficients $\alpha$ and $\delta$ only depends on $\lambda,b,r$ and $\eta$. 
Extension of this result to a family of Markov kernels $\{P_u:u\in [0,1]\}$ satisfying the previous conditions is then immediate. 
Then,  using the equivalence of the norms $\Vert\cdot\Vert_V$ and $\Vert\cdot\Vert_{V_{\delta}}$, one can show as in Proposition $2$ in \citet{Truquet2017} that there exists $C>0$ and $\rho\in (0,1)$ such that 
$$\sup_{u\in [0,1]}\Vert\delta_x P_u^n-\pi_u\Vert_V\leq CV(x)\rho^j.$$
Then the family of Markov kernels $\left\{P_u:u\in [0,1]\right\}$ is simultaneously $V-$uniformly ergodic.

\paragraph{Note.} 
The most important case for application of our results concerns the case $V_s=\phi^{q_s}$ for some $q_s\in (0,1)$ and $\phi:G\rightarrow [1,\infty)$ is a measurable function. See for instance Proposition \ref{sufficient} below for a result stated for this particular case. 
One can then obtain simultaneous $V_s-$uniform ergodicity for $s=0,1,\ldots,k$ if 
\begin{enumerate}
\item
there exists a positive real number $K$ such that for all $u\in [0,1]$, $P_u V_k\leq KV_k$,
\item
there exist an integer $m\geq 1$, two real numbers $\lambda\in (0,1)$, $b>0$, a family of positive real numbers $\{\eta_r:r>0\}$ and a family $\{\nu_r:r>0\}$ of probability measures on $G$ such that for all $r>0$ and all $u\in [0,1]$, the Markov kernel $P_u^m$ satisfies Condition $\mathcal{C}\left(V_k,\lambda,b,r,\eta_r,\nu_r\right)$. 
\end{enumerate}

Indeed, let $s\in \{0,1,\ldots,k-1\}$. From Jensen's inequality, we have, for any $s=0,\ldots,k$, $P_u V_s\leq K^{q_s/q_k}V_s$ and for any $r>0$, the family of Markov kernels $\left\{P^m_u: u\in [0,1]\right\}$ satisfies Condition $\mathcal{C}\left(V_s,\lambda^{q_s},b^{q_s},r^{q_s/q_k},\eta_r,\nu_r\right)$. From our previous discussion, we deduce that the family $\left\{P_u: u\in [0,1]\right\}$ is simultaneously $V_s-$uniformly ergodic.

\subsubsection{Other approach}
Simultaneous uniform ergodicity can also be obtained from other conditions. 
For instance, if each kernel $P_u$ is $V-$uniformly ergodic, then perturbation methods can be applied to get a local simultaneous
$V-$uniform ergodicity property which can easily be extended to the interval $[0,1]$ by compactness. When the Markov kernel is not continuous with respect to the operator norm, but satisfies some weaker continuity properties, \citet{Hervé} (see Theorem $1$) obtained a nice result based on the Keller-Liverani perturbation theorem.  The following proposition is an easy consequence of their result.  

\begin{prop}\label{central}
Let $\left\{P_u:u\in [0,1]\right\}$ be a family of Markov kernels on a measurable space $\left(G,\mathcal{G}\right)$ and $V:G\rightarrow [1,\infty)$ be a measurable function satisfying the three following conditions.

\begin{enumerate}
\item
For each $u\in[0,1]$, the Markov kernel $P_u$ admits a unique invariant measure $\mu_u$ such that $\int V d\mu_u <\infty$ and there exists $\kappa_u\in(0,1)$ and $C_u>0$ such that $\Vert \delta_x P_u^n-\mu_u\Vert_V\leq C_uV(x)\kappa_u^n$.

\item
There exist an integer $m\geq 1$, a real number $\lambda\in (0,1)$ and two positive real numbers $K$ and $L$ such that for all $u\in [0,1]$,
$$P_uV\leq K V,\quad P_u^mV\leq \lambda V+L.$$

\item
The application $u\rightarrow P_u$ is continuous for the norm $\Vert\cdot\Vert_{V,1}$.
\end{enumerate}

Then, there exists $\kappa\in (0,1)$ and $C>0$ such that 
$$\sup_{u\in [0,1]}\Vert \delta_x P_u^n-\mu_u\Vert_V=CV(x)\kappa^n.$$
Moreover $\sup_{u\in [0,1]}\Delta_V\left(P_u^n\right)=O\left(\kappa^n\right)$.
\end{prop}

\paragraph{Proof of Proposition \ref{central}}
Let $u\in [0,1]$. Our assumptions are exactly that of Theorem $1$ in \citet{Hervé}. This result guarantees the existence of an open interval $I_u\ni u$ of $[0,1]$, two real numbers $C_u>0$ and  $\kappa_u\in [0,1]$ such that for all $x\in G$, $\sup_{v\in I_u}\Vert \delta_x P_v^n-\mu_v\Vert_V\leq C_u V(x)\kappa_u^n$. From a compactness argument, $[0,1]$ can be covered by a finite number of such intervals $I_{u_1},\ldots,I_{u_p}$. Then the simultaneous $V-$uniform ergodicity condition follows by setting $\kappa=\max_{1\leq i\leq p}\kappa_{u_i}$ and defining the constant $C=\max_{1\leq i\leq p}C_{u_i}$. 
Moreover, we have (see (\ref{banddob})) 
$$\Delta_V\left(P_u^n\right)\leq \sup_{x\in G}\frac{\Vert \delta_x P_u^n-\mu_u\Vert_V}{V(x)},$$
which gives the second conclusion of the proposition.$\square$

\paragraph{Note.}
For the AR$(1)$ process $X_t(u)=\alpha(u)X_{t-1}(u)+\xi_t$, with $u\mapsto \alpha(u)\in (-1,1)$ continuous and $\xi_1$ has absolutely continuous error distribution with a density denoted by $\nu$ having a moment of order $a>0$, it is well known that  $P_u(x,dy)=\nu\left(y-\alpha(u)x\right)dy$ is $V_a-$geometrically ergodic with $V_a(x)=\left(1+\vert x\vert\right)^a$. See \citet{Hervé2}, Section $4$ for a discussion of the geometric ergodicity of some classical autoregressive processes. Moreover, the continuity of $u\mapsto P_u$ holds for the norm $\Vert\cdot\Vert_{V_a,1}$ as shown in \citet{Hervé}, Example $1$ (the result is shown for the case $a=1$ but extension to the case $a>0$ is straightforward). Then Proposition \ref{central} applies to this example.   
This approach does not require an additional property for the density $f_{\xi}$ such as existence of a positive lower bound on any bounded interval of the real line.  In contrast, positivity of the noise density is often required to check the small set condition in $\mathcal{C}\left(V,\lambda,b,R,\eta,\nu\right)$. However, construction of locally stationary Markov chain models considered in \citet{Truquet2017} is based on the simultaneous drift and small set conditions and we will not use Proposition \ref{central} in the rest of this paper.

\subsection{Regularity of higher-order finite dimensional distributions}
We now study existence of some derivatives for a functional $u\mapsto \int g d\pi_{u,j}$ where for $j\geq 2$, $g:E^j\rightarrow \R$ is a measurable function and $\pi_{u,j}(d{\bf x})=\pi_u(d x_1)Q_u(x_1,d x_2)\cdots Q_u(x_{j-1},dx_j)$. 
In time series analysis, a simple example is the estimation of the local covariance $u\mapsto \c\left(X_0(u),X_1(u)\right)$ where $\left(X_k(u)\right)_k$ is a stationary Markov chain with kernel $Q_u$.
For $x_1,\ldots,x_j\in E$ and $0\leq \ell\leq k$, we set $V_{\ell,j}(x_1,\ldots,x_j)=\sum_{i=1}^j V_{\ell}(x_i)$.
For an integer $j\geq 1$, we denote by $\mathcal{M}_V(E^j)$ the space of signed measures on $E^j$ such that 
$$\Vert \mu\Vert_V:=\sup\left\{\int f d\mu: \vert f(x_1,\ldots,x_j)\vert\leq V(x_1)+\cdots+V(x_j)\right\}.$$
Finally, let
$$M_{\ell}(x_1)=\sup_{u\in [0,1]}\int \left\vert \partial_1^{(\ell)}f(u,x_1,y_1)\right\vert \gamma(x_1,dy_1).$$

The following additional assumption will be needed.

\begin{description}
\item{\bf B4} There exists $C>0$ such that for $0\leq s\leq s+\ell\leq k$ and all $x_1,x_2\in E$, we have
$$V_s(x_1)M_{\ell}(x_2)\leq C\left(V_{s+\ell}(x_1)+V_{s+\ell}(x_2)\right).$$
\end{description}

Constant $C$ can be the same as in assumption {\bf B2}, this is why we use the same notation.
The following result is a consequence of Corollary \ref{sufficient}.

\begin{cor}\label{multiple}
Let $\left\{Q_u: u\in [0,1]\right\}$ be a family of Markov kernels on $E$ satisfying the assumptions {\bf B1-B4}. 
Then the application $u\mapsto \pi_{u,j}$ from $[0,1]$ to $\mathcal{M}_{V_0}(E^j)$, is $k-$times continuously differentiable.
\end{cor}

\paragraph{Note.} Assumption {\bf B4} will be satisfied if there exists a function $\phi:E\rightarrow [1,\infty)$
such that
$$\int \phi(y_j)^d\left\vert\partial^{(\ell)}_1f(u,x_j,y_j)\right\vert\gamma(x_j,dy_j)\leq C \phi(x_j)^{d+r_{\ell}}$$
for $0\leq d\leq d_0$. Indeed in this case, one can take (up to a constant) $V_{\ell}(x_j)=\phi(x_j)^{d+r\ell}$
with $r=\max\left(r_1,r_2/2,\ldots,r_k/k\right)$ and $k$ such that $d+rk\leq d_0$.

\paragraph{Proof of Corollary \ref{multiple}}
 Here we set for ${\bf x}\in E^j$ and $A\in\mathcal{E}^{\otimes j}$,
$$Q_{u,j}({\bf x},A)=\int_A f\left(u,x_j,y_j\right)\gamma_j({\bf x},d{\bf y}),$$
with $\gamma_j({\bf x},d{\bf y})=\gamma(x_j,dy_j)\prod_{i=1}^{j-1}\delta_{x_{i+1}}(dy_i)$. 
\begin{itemize}
\item
Let us first check that $Q_{u,j}$ satisfies assumption {\bf B1}. Let $1\leq s\leq k$.
For an integer $h\geq j$ and a measurable function $g:E^j\rightarrow \R$ such that $\vert g\vert\leq V_{s,j}$, we have 
\begin{eqnarray*}
\left\vert Q_{u,j}^jg({\bf x})\right\vert&\leq & \int \left\vert g(y_1,\ldots,y_j)\right\vert Q_u(x_j,dy_1)Q_u(y_1,dy_2)\cdots Q_u(y_{j-1},dy_j)\\
&\leq& \sum_{i=1}^j Q_u^iV_s(x_j)\\
&\leq& C_j V_s(x_j),
\end{eqnarray*}
with $C_j=\sum_{i=1}^j C^i$ and $C$ defined in (\ref{Hyp1}). 
We then get 
$$\sup_{\vert g\vert\leq V_{s,j}}\left\vert Q_{u,j}^hg({\bf x})-\pi_{u,j} g\right\vert\leq \sup_{\vert f\vert\leq C_j V_s}\left\vert Q_u^{h-j}f(x_j)-\pi_u f\right\vert\leq C_j\sup_{\vert f\vert\leq V_s}\left\vert Q_u^{h-j}f(x_j)-\pi_u f\right\vert.$$
From the simultaneous $V_s-$uniform ergodicity property for $\left\{Q_u:u\in [0,1]\right\}$, the previous bounds entail automatically {\bf B1}.

\item
Now assume that the family $\{Q_u:u\in [0,1]\}$ satisfies the assumptions {\bf B2-B3}. Then the family $\{Q_{u,j}:u\in [0,1]\}$
automatically satisfies the assumption {\bf B2} and {\bf B3} (\ref{Hyp2}). Let us check assumption {\bf B3} (\ref{Hyp1}). We have 
$$\int V_{s,j}({\bf y})\left\vert \partial_1^{(\ell)}f(u,x_j,y_j)\right\vert \gamma_j({\bf x},d{\bf y})\leq 
C\left[V_{s+\ell}(x_j)+\sum_{i=2}^jV_s(x_i)M_{\ell}(x_j)\right],$$
Using assumption {\bf B4}, we have $V_s(x_i)M_{\ell}(x_j)\leq C\left(V_{s+\ell}(x_i)+V_{s+\ell}(x_j)\right)$ and {\bf B3} (\ref{Hyp1}) is also satisfied for the family $\{Q_{u,j}:u\in [0,1]\}$. This completes the proof.$\square$
\end{itemize}

\section{Locally stationary Markov chains}
In this section, we consider a topological space $E$ endowed with its Borel $\sigma-$field $\mathcal{B}(E)$ and a triangular array of Markov chains $\left\{X_{n,t}:1\leq t\leq n, n\geq 1\right\}$ such that for all $(x,A)\in E\times \mathcal{B}(E)$ and $1\leq t\leq n$,   
$$\P\left(X_{n,t}\in A\vert X_{n,t-1}=x\right)=Q_{t/n}(x,A),\quad X_{n,0}\sim \pi_0.$$ 
We remind that for $u\in [0,1]$, $\pi_u$ denotes the invariant probability of $Q_u$.

\subsection{Some results about locally stationary Markov chains}
We first recall some results obtained in \citet{Truquet2017}. For simplicity, we introduce the two following conditions. For $\epsilon>0$, we denote $I_m(\epsilon)$ the subsets of $[0,1]^m$ such that $(u_1,\ldots,u_m)\in I_m(\epsilon)$ if and only if $\vert u_i-u_j\vert<\epsilon$ for $1\leq i,j\leq m$.
\begin{description}
\item[L1]
There exist a measurable function $V:E\rightarrow [1,\infty)$, an integer $m\geq 1$, some positive real numbers $\epsilon,K,\lambda,b,r,\eta$ with $\lambda<1$, $r>2b/(1-\lambda)$ and a probability measure $\nu$ such that for all $(u_1,u_2,\ldots,u_m)\in I_m(\epsilon)$, the kernel $Q_{u_1}Q_{u_2}\cdots Q_{u_m}$ satisfies Condition
$\mathcal{C}\left(V,\lambda,b,r,\eta,\nu\right)$. Moreover, there exists $K>0$ such that $Q_uV\leq KV$ for all $u\in [0,1]$.
\item[L2]
There exists a measurable function $V':E\rightarrow [1,\infty)$ such that $\sup_{u\in[0,1]}\pi_u V'<\infty$ and for all $x\in E$, $\Vert \delta_x Q_u-\delta_x Q_v\Vert_V\leq V'(x)\vert u-v\vert$. 
\item[L3]
For all $(u,v)\in [0,1]^2$, we have 
$$\Vert \delta_x Q_u-\delta_x Q_v\Vert_1\leq L(x)\vert u-v\vert,\mbox{ with } \sup_{u\in [0,1]\atop 1\leq \ell'\leq\ell}\E\left[ L\left(X_{\ell}(u)\right)V\left(X_{\ell'}(u)\right)\right]<\infty.$$
Here, $(X_t(u))_{t\in\Z}$ denotes a stationary time-homogeneous Markov chain with transition kernel $Q_u$.
\end{description} 
Under the conditions {\bf L1-L3}, it is shown in \citet{Truquet2017} (see Theorem $3$) that for all integer $j\geq 1$, the distribution $\pi_{t,j}^{(n)}$ of $\left(X_{n,t},\ldots,X_{n,t+j-1}\right)$ satisfies  
\begin{equation}\label{approximation}
\Vert \pi_{t,j}^{(n)}-\pi_{u,j}\Vert_V\leq C_j\left[\left\vert u-\frac{t}{n}\right\vert+\frac{1}{n}\right],
\end{equation}
where $C_j>0$ does not depend on $u,n,t$ and $V_j(x_1,\ldots,x_j)=V(x_1)+\cdots+V(x_j)$. 
Note that under Assumption {\bf L1} entails simultaneous $V-$uniform ergodicity for the family $\left\{Q_u:u\in [0,1]\right\}$ (See Section \ref{rat}).  
Condition {\bf L1} is useful to guarantee some $\beta-$mixing properties for the triangular array. See Proposition $3$ in \citet{Truquet2017} for details.
Note also that condition {\bf L2} is always satisfied for $(V,V')=(V_0,V_1)$ if assumption {\bf B3} (\ref{Hyp1}) holds true.
In \citet{Truquet2017}, Proposition $2$ and its proof, it has been shown that Assumptions {\bf L1-L3} entail, for each $u\in [0,1]$, geometric ergodicity of a Markov chain with transition kernel $Q_u$. Moreover, the finite dimensional distribution $\pi_{u,j}$ 
are shown to be Lipschitz with respect to $u$, when the space of signed measure on $E^j$ is endowed with the $V-$norm. However higher-order regularity (such as differentiability) has not been studied and this is precisely the aim of this paper.

For more clarity, we introduce the following terminology.
\begin{Def}\label{important}
A triangular array of Markov chains $\left\{X_{n,t}:1\leq t\leq n, n\geq 1\right\}$ associated to a family of Markov kernel $\{Q_u: u\in [0,1]\}$ will be said $V-$locally stationary if (\ref{approximation}) is satisfied.  
\end{Def}

\subsection{Simple sufficient conditions}
In order to check more easily our assumptions for specific examples, we give below a set of conditions that guarantee, for the same topology, local stationarity as well as differentiability of the applications $u\mapsto \pi_{u,j}$ for $j\geq 1$. In particular, the following set of assumptions will imply at the same time ${\bf L1-L3}$ and ${\bf B1-B3}$. Proposition \ref{sufficient} given below is then important for practical applications of our results to locally stationary Markov models.
We only consider the case of power functions, i.e. for each integer $s$, $V_s$ is a power of a measurable function $\phi:E\rightarrow [1,\infty)$. This is the most interesting case in practice. 

\begin{description}
\item[SC1]
There exist an integer $m\geq 1$, some positive real numbers $d_0,\epsilon,K,\lambda,b$ with $\lambda<1, d_0\geq 1$, a family of positive real number $\{\eta_r: r>0\}$ and a family $\{\nu_r: r>0\}$ of probability measures on $E$ such that for all $r>0$ and for all $(u_1,u_2,\ldots,u_m)\in I_m(\epsilon)$, the kernel $Q_{u_1}Q_{u_2}\cdots Q_{u_m}$ satisfies Condition
$\mathcal{C}\left(\phi^{d_0},\lambda,b,r,\eta_r,\nu_r\right)$. Moreover, there exists $K>0$ such that $Q_uV\leq KV$ for all $u\in [0,1]$.
\item[SC2]
There exists an integer $k\geq 1$ such that for all $(x,y)\in E^2$, the function $u\mapsto f(u,x,y)$ is $k-$times continuously differentiable.
\item[SC3]
There exist some real numbers $d_1>0$ and $q\geq 0$ such that $d_1+kq\leq d_0$ and for all $1\leq \ell\leq k$ and $d\leq d_1+(k-\ell)q$,  
$$\int \phi^d(y)\left\vert \partial_1^{(\ell)}f(u,x,y)\right\vert \gamma(x,dy)\leq C\phi^{d+q\ell}(x).$$
Moreover, for $s=0,\ldots,k$,
$$\lim_{h\rightarrow 0}\int \phi^{d_1+qs}(y)\left\vert \partial_1^{(k-s)}f(u+h,x,y)-\partial_1^{(k-s)}f(u,x,y)\right\vert
\gamma(x,dy)=0.$$
\end{description}

\begin{prop}\label{sufficient2}
Assume that {\bf SC1-SC3} hold true. Set $V_0=\phi^{d_1}$. 
The triangular array of Markov chain $\left\{X_{n,t}:1\leq t\leq n, n\geq 1\right\}$ is $V_0-$locally stationary. Moreover, for any integer $j\geq 1$, the application
$u\mapsto \pi_{u,j}$, from $[0,1]$ to $\mathcal{M}_{V_0}(E^j)$, is $k-$times continuously differentiable.
\end{prop}

\paragraph{Proof of Proposition \ref{sufficient2}}
For $s=0,\ldots,k$, we set $V_s=\phi^{d_1+qs}$. Note that from {\bf SC1}, Assumption {\bf L1} is automatically satisfied for each function $V_s$, $s=0,\ldots,k$.
Indeed if a Markov kernel $P$ satisfies for any $r>0$, the condition $\mathcal{C}(V,\lambda,b,r,\eta_r,\nu)$, then, for any $\kappa\in (0,1)$, it also satisfies condition $\mathcal{C}(V^{\kappa},\lambda^{\kappa},b^{\kappa},r^{\kappa},\eta_r,\nu)$.
 See the Note in Section \ref{explain} for a precise justification.

Moreover, from {\bf SC3} (set $d=d_1$ and $\ell=1$), Assumption {\bf L2} holds true for $V=V_0$ and $V'=V_1$. 

Next, we check {\bf L3}. Using {\bf SC3} with $d=0$ and $\ell=1$, we see that one can choose $L=C\phi^q$. Setting $V=V_0$, we know from {\bf L1-L2} that $\sup_{u\in [0,1]}\int \phi^{d_1+q}d\pi_u<\infty$. See \citet{Truquet2017}, Proposition $2$.
This shows that the integrability condition in {\bf L3} is satisfied. The proof of local stationary then follows.

We next check {\bf B1-B4}. {\bf B1} follows from {\bf L1} which holds true for all the functions $V_s$, $s=0,\ldots,k$. See the discussion of Section \ref{explain} for details.
Finally, {\bf B2-B4} follow directly from {\bf SC2-SC3}. See also the Note after Corollary \ref{multiple} for checking {\bf B4}. Differentiability of the marginal distributions then follows from Corollary \ref{multiple}. The proof is now complete.$\square$

\subsection{Application to bias control in nonparametric estimation}\label{locpol}
In this section, we discuss why differentiability properties of the application $u\mapsto \pi_{u,j}$ are fundamental for controlling the bias in nonparametric estimation of some parameter curves. 
Let $\left\{X_{n,t}:1\leq t\leq n, n\geq 1\right\}$ be a triangular array of $V-$locally stationary Markov chains. For a given integer $1\leq j\leq n$ and $1\leq t\leq n-j+1$, set $Z_{n,t}=\left(X_{n,t},\ldots,X_{n,t+j-1}\right)$. We also assume that the application $g:[0,1]\rightarrow \mathcal{M}_V(E^j)$ defined by $g(u)=\pi_{u,j}$ is $k-$times continuously differentiable. If Assumptions {\bf SC1-SC3} are satisfied, Proposition \ref{sufficient2}, given in the previous section, guarantees $V-$local stationarity et that $g$ is $k-$times continuously differentiable when $V=V_0$.

Let $f:E^j\rightarrow \R$ be a measurable function such that $\vert f\vert_V<\infty$. We want to estimate the quantity $\psi_f(u)=\int f d\pi_{u,j}$ using local polynomials. 
We precise that the approach used here is very classical in nonparametric estimation and, except for the local approximation, is identical to that used for i.i.d. data. See \citet{Tsy}, Section $1.8$, for a general approach for studying of the bias of local polynomial estimators.
Let $K$ be a continuous probability density, bounded and supported on $[-1,1]$ and $b\in(0,1)$ a bandwidth parameter such that $b=b_n\rightarrow 0$ and $nb\rightarrow \infty$.
We set $K_b=\frac{1}{b}K\left(\cdot/b\right)$.
An estimator $\hat{\psi}_f(u)$ of $\psi_f(u)$ is given by the first component of the vector 
$$\hat{\mathcal{H}}_f(u):=\left(\hat{\psi}_f(u),b\hat{\psi}_f'(u),\ldots,b^{k-1}\hat{\psi}_f^{(k-1)}(u)\right)'=\arg\min_{\alpha_0,\ldots,\alpha_{k-1}\in \R}\sum_{t=1}^nK_b\left(u-\frac{t}{n}\right)\left[f\left(Z_{n,t}\right)-\sum_{i=0}^{k-1}\alpha_i \frac{(t/n-u)^i}{b^i i!}\right]^2.$$
For $1\leq t\leq n$, we set 
$$v_t(u)=\left(1,\frac{t/n-u}{b},\ldots,\frac{(t/n-u)^{k-1}}{b^{k-1}(k-1)!}\right)'$$
and 
$$D(u)=\frac{1}{n-j+1}\sum_{t=1}^{n-j+1} K_b\left(t/n-u\right)v_t(u)v_t(u)',\quad  \hat{N}_f(u)=\frac{1}{n-j+1}\sum_{t=1}^{n-j+1} K_b\left(t/n-u\right)v_t(u)f\left(Z_{n,t}\right).$$
From (\ref{approximation}), we have 
$$\max_{1\leq t\leq n-j+1}\sup_{\vert f\vert_V\leq 1}\left\vert \E f\left(Z_{n,t}\right)-\psi_f(t/n)\right\vert=O\left(1/n\right).$$
Next, setting $\mathcal{H}_f(u)=\left(\psi_f(u),b\psi_f'(u),\ldots,b^{k-1}\psi_f^{(k-1)}(u)\right)'$ and using the differentiability properties of $\phi$, one can apply the bound (\ref{super}). There exists $C>0$ such that for all $n\geq 1$, $1\leq t\leq n-j+1$ and $u\in [0,1]$,  
$$\sup_{\vert f\vert_V\leq 1}\left\vert \psi_f(t/n)-\mathcal{H}_f(u)'v_t(u)\right\vert\leq C(u-t/n)^k.$$
We deduce that 
$$\sup_{u\in [0,1]}\sup_{\vert f\vert_V\leq 1}\left\vert \E\hat{N}_f(u)-D(u)\mathcal{H}_f(u)\right\vert=O\left(b^k+\frac{1}{n}\right).$$
The rest of the proof consists in bounding the matrix $D(u)^{-1}$ using very classical arguments available in the literature. Using our assumptions on the kernel and on the design $X_i=i/n$, the assumptions LP(1)-LP(3) of \citet{Tsy} are satisfied and Lemma $1.5$ and Lemma $1.7$ in \citet{Tsy} guaranty that 
 $\max_{u\in [0,1]}\Vert D(u)^{-1}\Vert=O(1)$. Then we get
$$\sup_{u\in [0,1]}\sup_{\vert f\vert_V\leq 1}\left\vert \E\hat{\mathcal{H}}_f(u)-\mathcal{H}_f(u)\right\vert=O\left(b^k+\frac{1}{n}\right).$$
In conclusion, up to a term of order $1/n$ which is negligible and can be interpreted as a deviation term with respect to stationarity, the bias is of order $b^k$ when $\psi_f$ is $k-$times continuously differentiable. We then recover a classical property of local polynomial estimators.

\paragraph{Notes}
\begin{enumerate}
\item
We will not discuss the variance of the estimator $\hat{\mathcal{H}}_f(u)$.
As shown in \citet{Truquet2017}, Proposition $3$, Assumption {\bf SC1} ensures geometric $\beta-$mixing properties for the triangular array of Markov chains $\left\{X_{n,t}: 1\leq t\leq n, n\geq 1\right\}$.   
Using standard arguments, one can then show that such variance is of order $1/nb$, as usual for nonparametric curve kernel estimators. Since this problem is not the scope of this paper, we omit the details.   
\item
Differentiability of $u\mapsto \pi_{u,2}$ is also important for deriving an expression of the bias for the local maximum likelihood estimator of some parameter curves. We defer the reader to Section $4.5$ in \citet{Truquet2017} for a discussion of this problem.
\end{enumerate}

\section{Extension to $p-$order Markov chains}
Let us now give an extension of our results to $p-$order Markov processes.
We choose here to present a version which can be applied directly to the examples of the last section of the paper. 
We consider a family $\left\{R_u: u\in [0,1]\right\}$ of probability kernel from $\left(E^p,\mathcal{E}^{\otimes p}\right)$
to $\left(E,\mathcal{E}\right)$. We assume that for $u\in [0,1]$, 
$$R_u\left({\bf x},A\right)=\int f\left(u,{\bf x},y\right)\gamma(dy),$$
for a measurable function $f:[0,1]\times E^{p+1}\rightarrow \R$ and a measure $\gamma$ on $E$.  
We also consider a triangular array $\left\{Y_{n,t}: 1\leq t\leq n, n\geq 1\right\}$ of $p-$order Markov processes such that
$$\P\left(Y_{n,t}\in A\vert Y_{n,t-1},\ldots,Y_{n,t-p}\right)=R_{t/n}\left(Y_{n,t-p},\ldots,Y_{n,t-1},A\right),\quad A\in \mathcal{B}(E),\quad 1\leq t\leq n.$$
For simplicity, we also define a sequence $\left(Y_{n,t}\right)_{t\leq 0}$ which is a time-homogeneous Markov process with transition kernel $R_0$.  
Note that setting $X_{n,t}=\left(Y_{n,t-p+1},\ldots,Y_{n,t}\right)$, one can define a triangular array $\left\{X_{n,t}: 1\leq t\leq n, n\geq 1\right\}$ of Markov chains. To this end, let $Q_u$ be the Markov kernel on $E^p$ defined by 
$$Q_u\left({\bf x},d{\bf y}\right)=R_u\left({\bf x},dy_p\right)\prod_{i=1}^{p-1}\delta_{x_{i+1}}(dy_i).$$
We then have 
$$\P\left(X_{n,t}\in A\vert X_{n,t-1}\right)=Q_{t/n}\left(X_{n,t-1},A\right),\quad A\in \mathcal{E}^{\otimes p},\quad 1\leq t\leq n.$$
One can then use the results available for locally stationary Markov chains to define and study some locally stationary Markov processes of order $p\geq 2$. For convenience, we give below a set of assumptions on the family of kernels $\{R_u:u\in [0,1]\}$ which ensure
local stationary and differentiability properties of $u\mapsto \pi_{u,j}$ for Markov chains with local transition kernels $\{Q_u:u\in [0,1]\}$. These properties will be derived using Proposition \ref{sufficient}.

Now let $\phi:E\rightarrow [1,\infty)$ a measurable function satisfying the following properties.
\begin{description}
\item {\bf SCp1} There exists a real number $d_0\geq 1$ and some positive real numbers $\alpha_{1,u},\ldots,\alpha_{p,u},\alpha_0$ such that $\sup_{u\in [0,1]}\sum_{i=1}^p\alpha_{i,u}<1$ and
$$R_u \phi^{d_0}\left({\bf x}\right)\leq \sum_{i=1}^p \alpha_{i,u} \phi^{d_0}(x_i)+\alpha_0,\quad {\bf x}\in E^p,\quad u\in [0,1].$$
Moreover, for each $r>0$, there exist a positive real number $\eta_r$ and a probability measure $\nu_r$ on $E$ such that,
$$R_u\left({\bf x}, A\right)\geq \eta_r \nu_r(A),\quad A\in \mathcal{E},\quad \max_{1\leq i\leq p}\phi(x_i)^{d_0}\leq r.$$

\item {\bf SCp2} There exists an integer $k\geq 1$ such that
for all $({\bf x},y)\in E^p\times E$, the function $u\mapsto f(u,{\bf x},y)$ is $k-$times continuously differentiable.

\item{\bf SCp3}
There exist some real numbers $d_1>0$ and $q\geq 0$ such that $d_1+ kq\leq d_0$ and for all $1\leq \ell\leq k$ and $d\leq d_1+(k-\ell)q$,  
$$\int \phi^d(y)\left\vert \partial_1^{(\ell)}f(u,{\bf x},y)\right\vert \gamma(dy)\leq C\sum_{i=1}^p \phi^{d+q\ell}(x_i).$$
Moreover,
$$\lim_{h\rightarrow 0}\int \phi^{d_1+qs}(y)\left\vert \partial_1^{(k-s)}f(u+h,{\bf x},y)-\partial_1^{(k-s)}f(u,{\bf x},y)\right\vert
\gamma(dy)=0.$$
\end{description}

\begin{cor}\label{pMarkov}
Assume that assumptions {\bf SCp1-SCp3} hold true. The triangular array of Markov chains $\left\{X_{n,k}:1\leq k\leq n, n\geq 1\right\}$ is $V_0-$locally stationary, with $V_0(x_1,\ldots,x_p)=\sum_{i=1}^p \phi^{d_1}(x_i)$. Moreover, for each integer $j\geq 1$, the finite dimensional distribution $u\mapsto \pi_{u,j}$ of the Markov chains with transition $Q_u$ it $k-$times continuously differentiable, as an application from $[0,1]$ to $\mathcal{M}_{V_0}(E^{pj})$.
\end{cor}

\paragraph{Proof of Corollary \ref{pMarkov}}
We will check the conditions of Proposition \ref{sufficient2}.
For an integer $j\geq 1$, and ${\bf x}\in E^j$, we set 
$$V_s({\bf x})=\sum_{i=1}^p \phi^{d_1+q s}(x_i),\quad 0\leq s\leq k.$$

\begin{enumerate}
\item
We first check the drift condition in {\bf SC1} for the function $V_k$. To this end, let ${\bf x}\in E^p$, $(u_i)_{i\geq 0}\in [0,1]^{\N}$ and $(Y_t)_{t\geq 1}$ a random sequence defined (on a given probability space) by 
$$\P\left(Y_n\in A\vert Y_{n-1},\ldots,Y_{n-p}\right)=R_{u_n}\left(Y_{n-p},\ldots, Y_{n-1},A\right),\quad n\geq p+1,\quad A\in\mathcal{E},$$  
and with arbitrary initial conditions $Y_i$, $1\leq i\leq p$.
It is then clear that the process $\left(X_n\right)_{n\geq p}$ defined by $X_n=\left(Y_n,\ldots,Y_{n-p+1}\right)$ for $n\geq p$
is a Markov chain of order $1$. We set $U_n=\E\left[\phi^{d_0}(Y_n)\vert X_p\right]$ for $n\geq 1$. 
From {\bf SCp1}, we have 
$$U_n\leq \sum_{i=1}^p \alpha_{i,u_n} U_{n-i}+\alpha_0,\quad n\geq p+1.$$
By induction, one can show that 
$$U_n\leq \alpha^{n/p}\max_{1\leq i\leq p}\phi^{d_0}(Y_i)+\frac{\alpha_0}{1-\alpha},$$
with $\alpha:=\sup_{u\in[0,1]}\sum_{i=1}^p \alpha_{i,u}$. This leads to the inequality
$$\E\left[V_k(X_n)\vert X_p={\bf x}\right]\leq \sum_{i=0}^{p-1}\alpha^{\frac{n-i+1}{p}}\max_{1\leq i\leq p}\phi^{d_0}(x_i)+\frac{p\alpha_0}{1-\alpha}.$$
Using the fact that $\alpha<1$, it is then clear that for an integer $m$ large enough, we have $\lambda:=\sum_{i=0}^{p-1}\alpha^{\frac{n-i+1}{p}}<1$ and then
$$Q_{u_{p+1}}\cdots Q_{u_{p+m}} V_k\leq \lambda V_k+ \frac{p\alpha_0}{1-\alpha}.$$
This shows the drift condition.

\item
Let us now check the small set condition for the function $V_k$. Let $r>0$ and assume that $V_k({\bf x})\leq r$. We then have $\max_{1\leq i\leq p}\phi^{d_0}(x_i)\leq r$. We set $\kappa_r=\nu_r\left(\left\{\phi^{d_0}\leq r\right\}\right)$. 
For $A\in \mathcal{E}^{\otimes p}$ and an integer $m\geq p$, we have using {\bf SCp1},
\begin{eqnarray*}
Q_{u_1}\cdots Q_{u_m}\left({\bf x},A\right)&=&\int \mathds{1}_A(x_{m+1},\ldots,x_{m+p})\prod_{i=p+1}^{m+p} R_{u_{i-p}}(x_{i-p+1},\ldots,x_{i-1},dx_i)\\
   &\geq& \int \mathds{1}_A(x_{m+1},\ldots,x_{m+p})\prod_{i=p+1}^{m+p} \mathds{1}_{\left\{\phi^{d_0}(x_i)\leq r\right\}}R_{u_{i-p}}(x_{i-p+1},\ldots,x_{i-1},dx_i)\\
	&\geq& \eta_r^m\int \mathds{1}_A(x_{m+1},\ldots,x_{m+p})\prod_{i=p+1}^{m+p} \mathds{1}_{\left\{\phi^{d_0}(x_i)\leq r\right\}}\nu_r(d x_i)\\
	&\geq& \eta_r^m\kappa_r^m\nu_{r,p}(A),
\end{eqnarray*}
with $\nu_{r,p}(A)=\kappa_r^{-p}\int_A \prod_{i=p+1}^{2p}\mathds{1}_{\left\{\phi^{d_0}(x_i)\leq r\right\}}\nu_r(dx_i)$.

\item
Next we check assumption {\bf SC3}. For $0\leq s\leq s+\ell\leq k$ and $x\in E^p$, we have from {\bf SCp3},
\begin{eqnarray*}
\int V_s(x_2,\ldots,x_p,y)\left\vert \partial^{(\ell)}_1 f(u,{\bf x},y)\right\vert \gamma(dy)
&\leq& C\left[\sum_{i=2}^p\phi^{d_1+qs}(x_i)\cdot\sum_{i=1}^p\phi^{q\ell}(x_i)+\sum_{i=1}^p\phi^{d_1+qs+q\ell}(x_i)\right]\\
&\leq& 3CV_{s+\ell}({\bf x}).
\end{eqnarray*}
In the last inequality, we used the bound $a^sb^{\ell}\leq a^{s+\ell}+b^{s+\ell}$ for $a,b\geq 1$.
The second part of {\bf SC3} follows directly from {\bf SCp3}. The result follows from Proposition \ref{sufficient2} and the proof is complete.$\square$
\end{enumerate}

\section{Examples}

In this section, we consider several examples of locally stationary Markov processes satisfying our assumptions and for which some parameter curves $u\mapsto \int f d\pi_{u,j}$ ($j\geq 1$)
can be estimated with local polynomials as explained in Section \ref{locpol}. 
We precise that our goal is not to estimate some parameter curves for the Markov kernel $Q_u=Q_{\theta(u)}$. 
However, as explained in Section \ref{locpol}, the results stated below are essential for getting an expression of the bias for minimum contrast estimators of $\theta(\cdot)$.
With respect to the examples discussed in \citet{Truquet2017}, Section \ref{logistique} provide a new example of locally stationary processes whereas Section \ref{aut} and Section \ref{INAR} give extensions to the order $p$ of some existing models.
A comparison of our results with that of \citet{Dahlhaus3} is given in Section \ref{aut}.

\subsection{Markov chains satisfying Doeblin's condition}
Here, we consider a family of Markov kernels $\left\{Q_u: u\in [0,1]\right\}$
such that for a probability measure $\mu$ and a measurable function $f:[0,1]\times E^2\rightarrow \R_+$, we have $Q_u(x,dy)=f(u,x,y)\mu(dy)$.

\begin{description}
\item[E11] 
There exists $c_{-}>0$ such that for all $(u,x,y)\in [0,1]\times E^2$, $f(u,x,y)\geq c_{-}$.
\item[E12]
There exists an integer $k\geq 1$ such that for all $(x,y)\in E^2$, the function $u\mapsto f(u,x,y)$ is of class $\mathcal{C}^k$ and 
$$\max_{0\leq \ell\leq k}\sup_{(u,x,y)\in [0,1]\times E^2}\left\vert \partial_1^{(\ell)}f(u,x,y)\right\vert<\infty.$$ 
\end{description} 

For a signed measure $\gamma$ on $E^j$, we define its total variation norm by 
$\Vert \gamma\Vert_1= \vert \gamma\vert (E)$. The total variation then coincides with the $V-$norm when $V\equiv 1$. 

\begin{prop}\label{Ex1}
Assume {\bf E11-E12}. Then the triangular array of Markov chain $\left\{X_{n,t}: 1\leq t\leq n\right\}$ is locally stationary for the total variation norm. Moreover, for all $j\geq 1$, 
the application $u\mapsto \pi_{u,j}$, as application from $[0,1]$ to $\mathcal{M}_1(E^j)$, is $k-$times continuously differentiable.
\end{prop}

\paragraph{Note.} This result is mainly interesting for compact state spaces $E$ (for instance with $E$ a compact subset of $\R^d$ and $\mu$ the uniform measure on $E$). Differentiability of the application $u\mapsto \pi_u$ can also be obtained from the results of \citet{Heid}. Assumption {\bf E11} entails Doeblin's condition, that is 
$Q_u(x,A)\geq c_{-}\mu(A)$ for all $(x,A)\in E\times\mathcal{B}(E)$.

\paragraph{Proof of Proposition \ref{Ex1}}
Local stationarity follows from \citet{Truquet2017} (see the second point in the Notes given after the statement of Theorem $1$). One can also use directly Proposition \ref{sufficient2}. {\bf SC1} is satisfied with $m=1$, $\phi\equiv 1$, $\eta_r=c_{-}$, $\nu_r=\mu$, $K=1$, $b=1$ and $\lambda=0$. Moreover, {\bf SC2-SC3} follows directly from {\bf E12}.$\square$

\subsection{Nonlinear autoregressive process}\label{aut}

We consider the following real-valued autoregressive process
$$X_{n,t}=m\left(t/n,X_{n,t-1},\ldots,X_{n,t-p}\right)+\sigma\left(t/n\right)\varepsilon_t,\quad 1\leq t\leq n,$$
where $m:[0,1]\times \R^p \rightarrow \R$ and $\sigma:[0,1]\rightarrow \R_+$ are two measurable functions and $\left(\varepsilon_t\right)_{t\in \Z}$ is a sequence of i.i.d random variables.
In what follows, we set $E=\R$ and for ${\bf y}\in \R^p$, $\vert {\bf y}\vert=\sum_{i=1}^p \vert y_i\vert$.
We will use the following assumptions.

\begin{description}
\item[E21]
The function $u\mapsto\sigma(u)$ is $k-$times continuously differentiable. Moreover $\sigma_{-}:=\inf_{u\in [0,1]}\sigma(u)>0$ and 
$$\max_{0\leq \ell\leq k}\sup_{u\in [0,1]}\left\vert \sigma^{(\ell)}(u)\right\vert<\infty.$$
\item[E22]
For all ${\bf y}\in \R^p$, the function $u\mapsto m(u,{\bf y})$ is $k-$times continuously differentiable. Moreover there exists a family of nonnegative real numbers $\left\{\beta_{i,u}: 1\leq i\leq p, u\in [0,1]\right\}$ such that $\sup_{u\in [0,1]}\sum_{i=1}^p\beta_{i,u}<1$ and four positive
real numbers $\beta_0,q',C_1,C_2$ such that for all $(u,{\bf y})\in [0,1]\times\R^p$,
$$\left\vert m(u,{\bf y})\right\vert\leq \sum_{i=1}^p\beta_{i,u}\vert y_i\vert+\beta_0,$$
$$\max_{1\leq \ell \leq k}\sup_{u\in [0,1]}\left\vert \partial_1^{(\ell)}m(u,{\bf y})\right\vert\leq C_1 \vert {\bf y}\vert^{q'}+C_2.$$
\item[E23]
The noise $\varepsilon_1$ has a moment of order $d_0$ such that $d_0-q'k>0$ and has a density $f_{\varepsilon}$, $k-$times continuously differentiable, positive everywhere and such that 
$$\int \vert y\vert^{d_0+(1-q')s}\left\vert f_{\varepsilon}^{(s)}(y)\right\vert dy<\infty,\quad s=0,\ldots,k.$$  
\end{description}
Setting 
$$R_u\left({\bf x},dy\right)=\frac{1}{\sigma(u)}f_{\varepsilon}\left(\frac{y-m(u,{\bf x})}{\sigma(u)}\right)dy,$$
the family $\left\{Y_{n,k}: 1\leq k\leq n, n\geq 1\right\}$ is a triangular array of time-inhomogeneous $p-$order Markov processes associated to the transition kernels $R_u$, $u\in [0,1]$.

\begin{prop}\label{autoreg}
Under the assumptions {\bf E21-E24}, the conclusions of Corollary \ref{pMarkov} hold true with $q=q'$, $d_1=d_0-q'k$ and $\phi(y)=1+\vert y\vert$, $y\in E$.
 \end{prop}

\paragraph{Example.} Consider the case for $p=1$ with $m(u,x)=\sum_{i=1}^I\left(a_i(u)x+b_i(u)\right)\mathds{1}_{x\in R_i}$, $\{R_1,\ldots,R_I\}$ a partition of $\R$ and $a_i,b_i$ are functions $k-$times continuously differentiable with $\max_{1\leq i\leq I}\max_{u\in [0,1]}\vert a_i(u)\vert<1$. 
This corresponds to a threshold model with non time-varying regions for the different regimes.
If {\bf E21} holds true, {\bf E22} follows with $q'=1$. If Assumption {\bf E23} is also valid for some $q_0>1$, Proposition \ref{autoreg} applies. This example is a generalization of the SETAR model discussed in \citet{Truquet2017} (see Example $3$ in Section $4.4$).

\paragraph{Proof of Proposition \ref{autoreg}}
To check the conditions of Corollary \ref{pMarkov}, we set $q=q'$ and
\begin{equation}\label{densities} 
f\left(u,{\bf x},y\right)=\frac{1}{\sigma(u)}f_{\varepsilon}\left(\frac{y-m(u,{\bf x})}{\sigma(u)}\right). 
\end{equation}
We first check the drift condition in {\bf SCp1}. 
Note first that
$$R_u\phi^{d_0}({\bf x})=\E\left[1+\left\vert m(u,{\bf x})+\sigma(u)\varepsilon_1\right\vert\right]^{d_0}.$$
From {\bf E22}, we then have 
$$R_u\phi^{d_0}({\bf x})\leq \E\left[\sum_{i=1}^p\beta_{i,u}\vert x_i\vert+1+\beta_0+\vert\varepsilon_1\vert\right]^{d_0}.$$
By convexity and setting $\beta=\sup_{u\in [0,1]}\sum_{i=1}^p\beta_{i,u}$, we get 
$$R_u\phi^{d_0}({\bf x})\leq \sum_{i=1}^p\beta_{i,u}\vert x_i\vert^{d_0}+\left(1-\beta\right)^{1-k}\E\left(1+\beta_0+\vert \varepsilon_1\vert\right)^{d_0}.$$
This shows the drift condition with 
$$\alpha_{i,u}=\beta_{i,u}\mbox{ and }\alpha_0=\left(1-\beta\right)^{1-k}\E\left(1+\beta_0+\vert \varepsilon_1\vert\right)^{d_0}.$$

Next, we check the small set condition. Suppose that ${\bf x}\in \R^p$ is such that $\phi(x_i)^{d_0}\leq r$, $1\leq i\leq p$. We then have 
$$\sup_{u\in [0,1]}\left\vert m(u,{\bf x})\right\vert\leq \sup_{u\in[0,1]}\sum_{i=1}^p\beta_{i,u}r^{1/{d_0}}+\beta_0.$$
Using the assumptions on $f_{\varepsilon}$ and $\sigma$, this entails that 
$$\widetilde{\eta}_r:=\inf_{\phi(x_i)\leq r, 1\leq i\leq p}\inf_{u\in[0,1]}\inf_{\vert y\vert\leq r}\frac{1}{\sigma(u)}f_{\varepsilon}\left(\frac{y-m(u,{\bf x})}{\sigma(u)}\right)>0.$$
We then get 
$$R_u({\bf x},A)\geq R_u({\bf x},A\cap [-r,r])\geq 2r\widetilde{\eta}_r\nu_r(A),$$
with $\nu_r$ the uniform distribution on $[-r,r]$. This shows the second part of {\bf SCp1}.

Assumptions {\bf E21-E23} and the expression (\ref{densities}) entail {\bf SCp2}. 

Let us now show {\bf SCp3}.
Using Assumptions {\bf E22-E23}, (\ref{densities}) and an induction argument, it can be shown that for $\ell=0,\ldots,k$,
\begin{equation}\label{shape}
\partial^{(\ell)}_1f\left(u,{\bf x},y\right)=\overline{\sigma}^{(\ell)}(u)f_{\varepsilon}\left(b_{y,{\bf x}}(u)/\sigma(u)\right)+\sum_{s=1}^{\ell}f_{\varepsilon}^{(s)}\left(b_{y,{\bf x}}(u)/\sigma(u)\right)\mathcal{P}_{\ell,u,s}\left(b_{y,{\bf x}}(u),b^{(1)}_{y,{\bf x}}(u),\ldots,b^{(\ell-s+1)}_{y,{\bf x}}(u)\right),
\end{equation}
with $\overline{\sigma}=1/\sigma$, $b_{y,{\bf x}}(u)=y-m(u,{\bf x})$ and for $1\leq s\leq \ell$, $\mathcal{P}_{\ell,u,s}$ is a polynomial of degree $s$ with coefficients of type $h(u)$ for bounded functions $h:[0,1]\rightarrow \R$.
One can then show that condition
$$\int \phi(y)^d \left\vert \partial_1^{(\ell)}f(u,{\bf x},y)\right\vert dy\leq C\sum_{i=1}^p\phi^{d+\ell q}(x_i)$$
holds true for $d\leq d_1+(k-\ell)q=d_0-\ell q$ if and only if $\int \vert y\vert^{d_0-q\ell+s}\cdot\left\vert f_{\varepsilon}^{(s)}(y)\right\vert dy<\infty$ when $s\leq \ell$. This is equivalent to 
$$\int \vert y\vert^{d_0+(1-q)s}\cdot\left\vert f_{\varepsilon}^{(s)}(y)\right\vert dy<\infty,\quad s=0,\ldots,k.$$
From {\bf E23}, the first part of {\bf SCp3} follows. To show the second part, we note that from the Lebesgue theorem, we have for each $M>0$,
$$\lim_{h\rightarrow 0}\int_{\{\vert y\vert\leq M\}}\phi^{d_1+qs}(y)\left\vert \partial_1^{(k-s)}f(u+h,{\bf x},y)-\partial_1^{(k-s)}f(u,{\bf x},y)\right\vert dy=0.$$ 
The second part of {\bf SCp3} will follow if we show that for ${\bf x}\in \R^p$,
\begin{equation}\label{montre}
\lim_{M\rightarrow \infty}\sup_{u\in [0,1]}\int_{\{\vert y\vert\geq M\}}\phi^{d_1+qs}(y)\left\vert \partial_1^{(k-s)}f(u,{\bf x},y)\right\vert dy=0.
\end{equation}
But one can show that (\ref{montre}) is a consequence of the expression of (\ref{shape}), the uniform integrability of $y\mapsto \phi(y)^{d_0+(1-q)s}f^{(s)}_{\varepsilon}(y)$ (which follows from {\bf E23})
and {\bf E21-E22}.

The result of the proposition is then a consequence of Corollary \ref{pMarkov}.$\square$

\paragraph{Notes}

\begin{enumerate}
\item
Let us compare our result with that of \citet{Dahlhaus3} who studied nonlinear autoregressive processes.
For simplicity, we restrict the study to $p=1$.
Suppose that for some $d_0\geq 1$, we have $\E\vert \varepsilon_1\vert^{d_0}<\infty$ and there exist $c>0$ and $\beta\in (0,1)$ such that
$$\sup_{u\in [0,1]}\vert m(u,x)-m(u,x')\vert\leq \beta\vert x-x'\vert,\quad \max_{i=1,2}\sup_{u\in[0,1}\left\vert \partial_i m(u,x)-\partial_im(u,x')\right\vert\leq C\vert x-x'\vert.$$
Theorem $4.8$ and Proposition $3.8$ in \citet{Dahlhaus3} show that the function $u\mapsto \int g d\pi_{u,j}$ is continuously differentiable whenever the function $g:E^j\rightarrow \R$ is continuously differentiable and satisfies for some $C>0$,
$$\vert g(z)-g(z')\vert\leq C(1+\vert z\vert^{d_0-1}+\vert z'\vert^{d_0-1})\vert z-z'\vert.$$
These authors also prove that there exists some positive constants $C_1$ and $C_2$ such that 
$$\E^{1/d_0}\vert X_{n,t}-X_t(u)\vert^{d_0}\leq C_1\left[\vert u-t/n\vert+1/n\right]\mbox{ with } X_t(u)=m\left(u,X_t(u)\right)+\sigma(u)\varepsilon_t,\quad t\in \Z,$$
and $\left\vert \int g d\pi^{(n)}_{t,j}-\int g d\pi_{u,j}\right\vert\leq C_2\left[\vert u-t/n\vert+1/n\right]$.
 
In contrast, when $d_0>1$, $\int \vert y\vert^{d_0}\left[f_{\varepsilon}(y)+\left\vert f_{\varepsilon}'(y)\right\vert\right]dy<\infty$ and there exist $\beta\in (0,1)$, $\beta',C>0$ such that 
$$\sup_{u\in [0,1]}\vert m(u,x)\vert\leq \beta\vert x\vert+\beta',\quad \sup_{u\in [0,1]}\vert \partial_1 m(u,x)\vert\leq C(1+\vert x\vert),$$ 
Proposition \ref{autoreg} guarantees that $u\mapsto \int gd\pi_{u,j}$ is continuously differentiable, provided that $\vert g(z)\vert\leq C(1+\vert z\vert)$ for some constant $C>0$.

One can then see that our assumptions on the regression function are less restrictive than that of \citet{Dahlhaus3} and no continuity assumption is made with respect to the second argument $x$. On the other hand, we impose much more regularity assumptions on the noise distribution (existence of a smooth density and a moment condition for its derivative). 
Our approach is interesting for non smooth regression functions. For instance, consider the threshold model with $m(u,x)=a_1(u)\max(x,0)+a_2(u)\max(-x,0)$ with $a_1,a_2$ continuously differentiable and $\max_{u\in [0,1]}\vert a_i(u)\vert<1$, $i=1,2$. In this case, the local approximation result of \citet{Dahlhaus3} is still valid but differentiability of $u\mapsto \int g d\pi_{u,j}$ cannot be obtained. In contrast, we can prove this differentiability for a different class of functions. More general threshold models with a discontinuous regression function can also be considered as in the example given just after the statement of Proposition \ref{autoreg}.

When the assumptions of \citet{Dahlhaus} are satisfied, our method also provides approximation of $\int g d\pi^{(n)}_{t,j}$ and smoothness of $u\mapsto \int g d\pi_{u,j}$ for very irregular functions (for instance indicators of Borel sets).

To conclude, we see on this particular example that for autoregressive processes, our results can afford a complement to that of \citet{Dahlhaus3}, for studying local approximation and smoothness properties of non smooth functions $g$ or for studying non smooth regression functions. Note also that we also provide a criterion for higher-order differentiability, a problem not considered in \citet{Dahlhaus3}.
\item

Exponential stability can be used for such models provided that $f_{\varepsilon}^{(s)}$ has some exponential moments for $s=0,\ldots,k$. In this case, one can take $V_s(y)=\exp\left(\kappa \vert y\vert\right)$ for all $s$. A precise result is given in Proposition \ref{expo} given in the appendix. In this case, the approach of \citet{Heid} can also be used for studying existence of derivatives and our general result, which also covers this case, is not useful (except that we provide a criterion for $p-$order Markov chain, which is new). However, these exponential moments induce a serious restriction on the noise distribution because fatter tails distributions such as Student distributions are excluded.
However, the local stationarity property of this model, resulting from Proposition \ref{expo}, is a new result.

\item
Our result can be also applied to the AR$(1)$ process $X_t= \alpha X_{t-1}+\varepsilon_t$ for getting derivatives of the applications $\alpha \mapsto \pi_{\alpha}$, as in \citet{Hervé}.
The index $u$ is replaced with $\alpha$ and the interval $[0,1]$ with $I=[-1+\epsilon,1-\epsilon]$ for some $\epsilon\in (0,1)$. Let $Q_{\alpha}(x,dy)=f_{\varepsilon}(y-\alpha x)dy$.
 In this case, one can take $q'=1$, $\phi(x)=1+\vert x\vert$ and if $k<d_0<k+1$, $d_1=d_0-k$. 
Under some assumptions that guaranty {\bf E23}, \citet{Hervé} showed in their Proposition $1$ that $\alpha\mapsto \pi_{\alpha}$, considered as an application from $I$ to $\mathcal{\phi^{\beta}}(\R)$,
is $k-$times continuously differentiable, provided that $0<\beta<d_1$. See their condition on $\beta$ given after the statement of their Lemma $1$. One can then see that our result is stronger. We claim that the slight difference between the two results is explained by the additional topologies used in their Lemma $1$ for studying continuity of the application 
$$\alpha\mapsto Q^{(\ell)}_{\alpha}(x,dy)=(-1)^{\ell}x^{\ell}f_{\varepsilon}^{(\ell)}(y-\alpha x)dy.$$
Let us enlighten why by supposing that $k=1$. From Theorem \ref{general}, we have, using our notations $T_{\alpha}^{(\ell)}\mu=\mu Q_{\alpha}^{(\ell)}$,
$$\pi_{\alpha}^{(1)}=(I-T_{\alpha})^{-1}T_{\alpha}^{(1)}\pi_{\alpha}.$$
Denoting by $\mathcal{L}\left(\phi^{d_1},\phi^{d_1'}\right)$ the set of bounded linear operators from $\mathcal{M}_{\phi^{d_1}}(\R)$ to $\mathcal{M}_{\phi^{d_1'}}(\R)$, the application $\alpha\mapsto T_{\alpha}^{(1)}$, as an application from $I$ to $\mathcal{L}\left(\phi^{d_1},\phi^{d_1'}\right)$ is only continuous when $d_1'<d_1$. This shows that one can only get continuity  $\alpha\mapsto \pi^{(1)}_{\alpha}$ for $\Vert\cdot\Vert_{\phi^{d_1'}}$ if we use operator norms. 
On the other hand, if $\mu\in\mathcal{M}_{\phi^{d_1}}(\R)$, one can show that the application $\alpha\mapsto T^{(1)}_{\alpha}\mu$, as an application from $I$ to $\mathcal{M}_{\phi^{d_1}}(\R)$ is continuous. As shown in Theorem \ref{general}, this weaker continuity condition is sufficient for getting continuity of $\alpha\mapsto \pi^{(1)}_{\alpha}$, as an application from $I$ to $\mathcal{M}_{\phi^{d_1}}(\R)$.   
\end{enumerate}

\subsection{Integer-valued time series}\label{INAR}
For $u\in [0,1]$ and $1\leq i\leq p$, let $\zeta_{i,u}$ and $\xi_u$ be some probability distributions supported on the nonnegative integers and 
for ${\bf x}\in\Z_+^p$, $R_u({\bf x},\cdot)$ will denote the probability distribution given by the convolution product
$\zeta_{1,u}^{*x_1}*\zeta_{2,u}^{*x_2}*\cdots *\zeta_{p,u}^{*x_p}*\xi_u$ with $\zeta_{i,u}^{*x}=\zeta_{i,u}^{*(x-1)}*\zeta_{i,u}$ if $x\geq 1$, $\zeta_{i,u}^{*1}=\zeta_{i,u}$ and the convention $\zeta_{i,u}^{*0}=\delta_0$.

Let us comment this Markov structure. When $p=1$, $R_u$ is the transition matrix of a Galton-Watson process with immigration.
Such Markov processes are also used in time series analysis of discrete data. 
For instance, if $\zeta_{i,u}$ denotes the Bernoulli distribution of parameter $\alpha_{i,u}$, 
such Markov processes are called INAR processes and were studied in \citet{Al} and \citet{Guan}. 
Note that in this case, we have the autoregressive representation 
$X_k=\sum_{i=1}^p\alpha_{i,u}\circ X_{k-i}+\varepsilon_k$, where $\alpha\circ x$ denotes a random variable following a binomial distribution of parameters $(x,\alpha)$ and independent from $\varepsilon_k$, an integer-valued random variable with probability distribution $q_u$. When $\zeta_{i,u}$ denotes the Poisson distribution of parameter $\alpha_{i,u}$ and $\xi_u$ denotes the Poisson distribution of parameter $\alpha_{0,u}$, then $R_u({\bf x},\cdot)$ is the Poisson distribution of parameter $\alpha_{0,u}+\sum_{i=1}^p\alpha_{i,u}$ and the Markov process coincides with the INARCH process studied in \citet{Ferland}.    
The distributions $\zeta_{i,u}$ and $\xi_u$ can also have a general form as in the generalized INAR processes studied by \citet{Latour} and are not
required to have exponential moments. For instance the log-logistic distribution $\zeta$ with parameters $\alpha,\beta>0$ and defined by $\zeta(x)=(\beta/\alpha)(x/\alpha)^{\beta}\left(1+(x/\alpha)^{\beta}\right)^{-2}$ for $x\in\Z_+$, has only a finite moment of order $k<\beta$.  
When $p=1$, conditions ensuring local stationarity for the INARCH and INAR processes are discussed in \citet{Truquet2017}.
Here, we propose an extension to the case $p\geq 1$, with general probability distributions $\zeta_{i,u}$and $\xi_u$ and additionally, we study the regularity properties of the marginal distributions w.r.t. $u$, a problem which has not been addressed before. 

We will use the following assumptions.

\begin{description}

\item[E31]
We have $\alpha:=\sup_{u\in [0,1]}\sum_{i=1}^p\sum_{x\geq 0}x \zeta_{i,u}(x)<1$ and there exists an integer $x_0$ such that $\beta:=\inf_{u\in[0,1]}\xi_u(x_0)>0$.

\item[E32]
For each integer $x\geq 0$, the applications $u\mapsto \zeta_{i,u}(x)$ and $u\mapsto \xi_u(x)$ are of class $\mathcal{C}^k$. Moreover, there exists a positive integer $d_1$ such that for $s=0,1,\ldots,k$,
$$\lim_{M\rightarrow \infty}\sup_{u\in[0,1]}\sum_{i=1}^p\sum_{x\geq M}x^{d_1+k-s}\left[\vert \zeta_{i,u}^{(s)}(x)\vert+\vert \xi_u^{(s)}(x)\vert\right]=0.$$

\end{description}

\begin{prop}\label{GW}
Assume that the assumptions {\bf E31-E32} hold true and set $\phi(x)=1+x$ for $x\in\N$ and $d_0=d_1+k$.
Then the conclusions of Corollary \ref{pMarkov} hold true.
\end{prop}

\paragraph{Note.} Assumption ${\bf E32}$ is satisfied for Bernoulli, Poisson or negative binomial distributions provided the real-valued parameter of these distributions is a $\mathcal{C}^k$ function taking values in the usual intervals $(0,1)$ (for the Bernoulli or negative binomial distribution) or $(0,\infty)$ (for the Poisson distribution).

\paragraph{Proof of Proposition \ref{GW}}
We check Assumptions {\bf SCp1-SCp3} of Corollary \ref{pMarkov}.

\begin{itemize}
\item
Let us first check the small set condition in {\bf SCp1}.  This small set condition is in fact satisfied for any finite set $C=\{0,1,\ldots,c\}^p$. Indeed, from assumption {\bf E31}, we have 
$\zeta_{i,u}(0)\geq 1-\alpha$ for $1\leq i\leq p$ and for ${\bf x}\in C$ and $A\in\mathcal{P}(\Z_+)$,
$$R_u({\bf x},A)\geq \left(1-\alpha\right)^{cp} \beta \delta_{x_0}(A).$$

\item
Next, we check the drift condition in {\bf SCp1} for the function $\phi^{d_0}$.
We denote by $\Vert\cdot\Vert_{d_0}$ the standard norm for the space $\L^{d_0}$. We remind that for some independent random variables $Y_1,\ldots,Y_n$ with mean $0$, the Burkhölder inequality gives the bound
$$\Vert\sum_{i=1}^n Y_i\Vert_{d_0}\leq C_{d_0}\Vert \sqrt{\sum_{i=1}^n Y_i^2}\Vert_{d_0},$$
where $C_{d_0}>0$ only depends on $d_0$. This leads to the bound
$$\Vert\sum_{i=1}^n Y_i\Vert_{d_0}\leq \left\{\begin{array}{c} C_{d_0}\sqrt{\sum_{i=1}^n \Vert Y_i\Vert_{q_0}^2}\mbox{ if } q_0\geq 2,\\ C_{d_0}\left(\sum_{i=1}^n \Vert Y_i\Vert_{d_0}^{d_0}\right)^{1/d_0}\mbox{ if } q_0\in (1,2)\end{array}\right.$$
If $\max_{1\leq i\leq n}\Vert Y_i\Vert_{d_0}\leq \kappa$, we then obtain
\begin{equation}\label{bound}
\Vert \sum_{i=1}^n Y_i\Vert_{d_0}\leq C_{d_0}\kappa n^{1/d_0'} \mbox{ with } d'_0:=\max(d_0,2).
\end{equation}
Next we set 
$$\kappa^{d_0}=\sup_{u\in [0,1]}\left[\max_{1\leq i\leq p}\sum_{x\geq 0}x^{d_0}\zeta_{i,u}(x)\right]\vee \sum_{x\geq 0}x^{d_0}\xi_u(x),$$
which is finite from assumption {\bf G2}. 
Let $S_u({\bf x})$ be a random variable with distribution $R_u\left({\bf x},\cdot\right)$. Since $S_u({\bf x})$ can be represented as a sum of $n=1+\sum_{i=1}^p x_i$ independent random variables, we deduce that   
$$\Vert S_u({\bf x})-\E S_u({\bf x})\Vert_{d_0}\leq C_{d_0}\kappa \left(1+\sum_{i=1}^px_i\right)^{1/d_0'}.$$
Setting $m_{i,u}=\sum_{x\geq 0} x \zeta_{i,u}(x)$, we have $\E S_u({\bf x})=\sum_{i=1}^p \zeta_{i,u}x_i$ and we get
$$\Vert 1+S_u({\bf x})\Vert_{d_0}\leq 2+\sum_{i=1}^p m_{i,u}x_i+ C_{d_0}\kappa \sum_{i=1}^p x_i^{1/d_0'}.$$
For any $\varepsilon>0$, there exists $b>0$ (depending on $\varepsilon$ and $d_0$) such that  $ C_{d_0}x_i^{1/d_0'}\leq \varepsilon x_i+b$.
We choose $\varepsilon$ such that $\alpha+2\varepsilon<1$. We then obtain, setting $b'=2+pb$ and $m_{i,u,\varepsilon}=m_{i,u}+\varepsilon$,
$$\Vert 1+S_u({\bf x})\Vert_{d_0}\leq \sum_{i=1}^p m_{i,u,\varepsilon}x_i+b'.$$
Using the equality $b'=\varepsilon b'/\varepsilon$ and convexity, we then obtain
$$R_u\phi^{d_0}({\bf x})=\Vert 1+S_u({\bf x})\Vert^{d_0}_{d_0}\leq \left(\alpha+2\varepsilon\right)^{d_0-1}\left[\sum_{i=1}^pm_{i,u,\varepsilon}x_i^{d_0}+\varepsilon^{1-k}(b')^{d_0}\right].$$
We the deduce that the drift condition in {\bf SCp1} is satisfied with 
$$\alpha_0=\left(\alpha+2\varepsilon\right)^{d_0-1}\varepsilon^{1-k}(b')^{d_0},\quad \alpha_{i,u}=\left(\alpha+2\varepsilon\right)^{d_0-1}m_{i,u,\varepsilon}.$$

\item 
Next we check the first part of {\bf SCp3}. 
In the sequel we set for an integer $M\geq 0$,
$$\mathcal{D}_M:=\sup_{u\in[0,1]}\max_{1\leq i\leq p}\max_{s=0,1,\ldots,k}\sum_{x\geq M}(1+x)^{d_1+k-s}\left[\vert \zeta_{i,u}^{(s)}(x)\vert+\vert \xi_u^{(s)}(x)\vert\right].$$
For ${\bf x}\in \Z_+^p$, the conditional density $f(u,{\bf x},\cdot)$ with respect to the counting measure on $\Z_+$ is given by the convolution product $\zeta_{1,u}^{*x_1}*\zeta_{2,u}^{*x_2}*\cdots *\zeta_{p,u}^{*x_p}*\xi_u$. Setting $n=x_1+\cdots+x_p+1$, we have for $1\leq \ell\leq k$,
$$\partial_1^{(\ell)}f(u,{\bf x},y)=\sum_{j_1+\cdots+j_n=y}\sum_{\ell_1+\ell_2+\cdots+\ell_n=\ell}\frac{\ell!}{\prod_{i=1}^n\ell_i!}p^{(\ell_1)}_{1,u}(j_1)\cdots p_{n,u}^{(\ell_n)}(j_n),$$
where $p_{1,u}=\cdots=p_{x_1,u}=\zeta_{1,u}$, $p_{x_1+1,u}=\cdots=p_{x_1+x_2,u}=\zeta_{2,u}$ and so on, up to
$p_{n-x_p+1,u}=\cdots=p_{n,u}=\xi_u$.
If $0\leq s\leq s+\ell\leq k$, we have, using the convexity of the application $y\mapsto V_s(y):=1+y^{d_1+s}$,
$$\sum_{y\geq 0}V_s(y)\left\vert\partial_1^{(\ell)}f(u,{\bf x},y)\right\vert\leq
\sum_{\ell_1+\ell_2+\cdots+\ell_n=\ell}\frac{\ell!}{\prod_{i=1}^n\ell_i!}n^{d_1+s}\mathcal{D}_0^{\ell+1}=n^{d_1+s+\ell}\mathcal{D}_0^{\ell+1}.$$
In the previous inequality, we have used the following property: if $\ell_1+\ell_2+\cdots+\ell_n=\ell$ then at most $\ell$ of these integers are positive and then for $i=1,2,\ldots,n$, 
$$\sum_{j_1,\ldots,j_n\geq 0}V_s(j_i)\left\vert p^{(\ell_1)}_{1,u}(j_1)\cdots p^{(\ell_n)}_{u,n}(j_n)\right\vert
\leq \mathcal{D}_0^{\ell+1}.$$
Using again convexity, we have $n^{d_1+\ell+s}\le p^{d_1+\ell+s-1}\sum_{i=1}^pV_{s+\ell}(x_i)$ and the first part of {\bf SCp3} follows.

\item
Finally, we check the second part of {\bf SCp3}. Let $x$ be a nonnegative integer and $s$ an integer such that $0\leq s\leq k$. It is easily seen that 
$$\lim_{h\rightarrow 0}\sum_{y\leq M}\left\vert\partial_1^{(k-s)}f(u+h,{\bf x},y)-\partial_1^{(k-s)}f(u,{\bf x},y)\right\vert V_s(y)=0.$$
Then it remains to show that
\begin{equation}\label{reste}
\lim_{M\rightarrow \infty}\sup_{u\in [0,1]}\sum_{y\geq M}\left\vert\partial_1^{(k-s)}f(u,{\bf x},y)\right\vert V_s(y)=0.
\end{equation}
But as for the proof of the first part of {\bf SCp3}, we have
$$\sum_{y\geq M}V_s(y)\left\vert\partial_1^{(k-s)}f(u,{\bf x},y)\right\vert\leq 
\sum_{\ell_1+\ell_2+\cdots+\ell_n=k-s}\frac{(k-s)!}{\prod_{i=1}^n\ell_i!}n^{d_1+s+1}\mathcal{D}_0^{k-s}\mathcal{D}_{M/n}=n^{d_1+k+1}\mathcal{D}_0^{k-s}\mathcal{D}_{M/n}.$$
The single change is the bound
\begin{eqnarray*}
\sum_{j_1+\cdots+j_n\geq M}V_s(j_i)\left\vert p^{(\ell_1)}_{1,u}(j_1)\cdots p^{(\ell_n)}_{u,n}(j_n)\right\vert
&\leq& \sum_{a=1}^n\sum_{j_1,\ldots,j_n\geq 0}\mathds{1}_{j_a\geq M/n}V_s(j_i)\left\vert p^{(\ell_1)}_{1,u}(j_1)\cdots p^{(\ell_n)}_{u,n}(j_n)\right\vert\\
&\leq& n\mathcal{D}_0^{k-s}\mathcal{D}_{M/n}.
\end{eqnarray*}
From the assumption {\bf E32}, we get (\ref{reste}), which completes the proof.$\square$

\end{itemize}

\subsection{Markov chain in a Markovian random environment}\label{logistique}
We consider a state space $E=E_1\times E_2$ with $E_1$ a finite set and $E_2$ an arbitrary metric space.
Let $\left\{P(u,\cdot,\cdot;z): u\in [0,1], z\in E_2\right\}$ is a family of stochastic matrices on $E_1$ and
$\left\{\overline{Q}_u:u\in [0,1]\right\}$ a family of Markov kernels on $E_2$. 
We assume that for all $u\in [0,1]$, $\overline{Q}_u(x_2,dy_2)=\overline{f}(u,x_2,y_2)\overline{\gamma}(x_2,dy_2)$ for 
a measurable function $\overline{f}:[0,1]\times E_2^2\rightarrow \R_+$ and a measure kernel $\overline{\gamma}$
on $E_2$. 
We consider the family of Markov kernels $\left\{Q_u:u\in [0,1]\right\}$ such that
$$Q_u((y_1,z_1),(d y_2,dz_2))=P\left(u,y_1,y_2;z_2\right)\overline{Q}_u(z_1,dz_2),\quad u\in [0,1].$$
Setting $f(u,(y_1,z_1),(y_2,z_2))=P(u,y_1,y_2;z_2)\overline{f}(u,z_1,z_2)$ we have 
$$Q_u((y_1,z_1),(d y_2,dz_2))=f(u,(y_1,z_1),(y_2,z_2))c(dy_2)\overline{\gamma}(z_1,dz_2),$$
where $c$ denotes the counting measure on $E_1$. We then set $\gamma\left((y_1,z_1),(dy_2,dz_2)\right)=c(dy_2)\overline{\gamma}(z_1,dz_2)$.

For $u\in [0,1]$, $P(u,\cdot,\cdot;z)$ is the transition matrix of a process in a Markovian random environment.
The kernels $Q_u$ can also be seen as a transition operator for a categorical time series with exogenous covariates. Indeed, if $\left\{X_{n,t}=\left(Y_{n,t},Z_{n,t}\right): 1\leq t\leq n, n\geq 1\right\}$ is a triangular array associated to the family $\left\{Q_u: u\in [0,1]\right\}$, we have 
$$\P\left(Y_{n,t}=y'\vert Y_{n,t-1},Z_{n,1},\ldots,Z_{n,n}\right)=P\left(t/n,y,y';Z_{n,t}\right),\quad 1\leq t\leq n, n\geq 1.$$
In the time-homogeneous case, \citet{FT1} recently studied Markov chains models with exogenous covariates of a general form and discussed their link with Markov chains in a random environment for studying ergodicity properties. We provide here a locally stationary analogue but with a restriction on the covariate process which is given by a locally stationary Markov chain. An important important example of such models is the autoregressive logistic model with $E_1=\{0,1\}$, $E_2=\R^g$ and 
$$P(u,y,1,z)=\frac{\exp\left(a_0(u)+a_1(u)y+z'\beta(u)\right)}{1+\exp\left(a_0(u)+a_1(u)y+z'\beta(u)\right)}$$
for some continuous functions $a_0,a_1:[0,1]\rightarrow \R$ and $\beta:[0,1]\rightarrow \R^g$.

We will use the following set of assumptions.
\begin{description}
\item[E41] For all $(y_1,y_2,z_2)\in E_1^2\times E_2$, the functions $u\mapsto P_u(y_1,y_2;z_2)$ is $k-$times continuously differentiable and positive. 

\item[E42] The family of Markov kernels $\left\{\overline{Q}_u: u\in [0,1]\right\}$ satisfies Assumptions {\bf SC1-SC3}. We denote the different constants involved in the assumptions by an overline.

\item[E43] For $\ell=0,\ldots,k$, we have 
$$\sup_{u\in [0,1]}\max_{y,y'\in E_1}\left\vert\partial_1^{(\ell)}P(u,y,y';z)\right\vert\leq C\overline{\phi}(z)^{\overline{q}\ell}.$$

\end{description}

We set $\phi(y,z)=\overline{\phi}(z)$.

\begin{prop}\label{randomenv}
Under the assumptions {\bf E41-E43}, the conclusions of Proposition \ref{sufficient} are valid for $(d_0,d_1,q)=\left(\overline{d}_0,\overline{d}_1,\overline{q}\right)$.

\end{prop}

\paragraph{Proof of Proposition \ref{randomenv}}
We check the condition of Proposition \ref{sufficient2}. 
To this end, we set $m=\overline{m}$ and $\epsilon=\overline{\epsilon}$.
The drift condition of {\bf SC1} is satisfied using our assumptions on the family $\left\{\overline{Q}_u: u\in [0,1]\right\}$. Indeed, we have $Q_{u_1}\cdots Q_{u_m}\phi^{d_0}=\overline{Q}_{u_1}\cdots\overline{Q}_{u_m}\phi^{d_0}$.

Next, let us check the small set condition in {\bf SC1}. We set $g(z)=\min_{y,y'\in E_1}\inf_{u\in [0,1]}P_u(y,y';z)$. From Assumption {\bf E41}, the function $g$ is positive. Let $h:E\rightarrow [0,1]$ a measurable function. If $\phi^{d_0}(y_0,z_0)\leq r$ for some $r>0$, we have 
\begin{eqnarray*}
Q_{u_1}\cdots Q_{u_m}h(y_0,z_0)&=&\int h(y_m,z_m)\prod_{i=1}^m P_{u_{m-i+1}}(y_{i-1},dy_i;z_i)\overline{Q}_{u_{m-i+1}}(z_{i-1},dz_i)\\
&\geq& \sum_{y_1\in E_1}\int h(y_m,z_m)g(z_m)\overline{Q}_{u_1}\cdots \overline{Q}_{u_m}(z_0,d z_m)\\
&\geq &\overline{\eta}_r\sum_{y_1\in E_1}\int h(y_m,z_m)g(z_m)\overline{\nu}_r(dz_m).
\end{eqnarray*}
This shows the drift condition with $\nu_r(dy,dz)=\vert E_1\vert^{-1} g(z)\overline{\nu}_r(dz)c(dy)/\int g(z')\overline{\nu}(d z')$, $c$ the counting measure on $E_1$ and $\eta_r=\vert E_1\vert\overline{\eta}_r\int g(z')\overline{\nu}(dz')$.   

Finally, we check {\bf SC3}. For $\ell=0,\ldots,k$, $0\leq \ell'\leq \ell$ and $d\leq d_1+(k-\ell)$, Assumptions {\bf E42-E43} guarantee that for a suitable positive constant $C$,
$$\int \overline{\phi}(z')^d\left\vert \partial_1^{(\ell')}P(u,y,y';z')\partial_1^{(\ell-\ell')} 
\overline{f}(u,z,z')\right\vert c(dy')\overline{\gamma}(z,dz')\leq C\phi(z)^{d+q\ell}.$$
Using the general Leibniz rule for the derivative of a product of functions, we get 
the first part of {\bf SC3}. To get the second part, let
$$A_h=\int \overline{\phi}^{d_1+qs}(z')\left\vert \partial_1^{(\ell')}P(u+h,y,y';z')-\partial_1^{(\ell')}
P(u,y,y';z')\right\vert\cdot\left\vert \partial_1^{k-s-\ell'}\overline{f}(u,z,z')\right\vert c(dy')\overline{\gamma}(z,dz')$$
and
$$B_h=\int \overline{\phi}^{d_1+qs}(z')\left\vert \partial_1^{(k-s-\ell')}\overline{f}(u+h,z,z')-\partial_1^{(k-s-\ell')}
\overline{f}(u,z,z')\right\vert\cdot\left\vert \partial_1^{k-s-\ell'}P(u+h,y,y';z')\right\vert c(dy')\overline{\gamma}(z,dz').$$
Using {\bf E43} and {\bf SC3} for $\overline{f}$, we have $\lim_{h\rightarrow 0}B_h=0$. 
Moreover, from {\bf E41}, {\bf E43}, {\bf SC3} for $\overline{f}$ and the Lebesgue theorem, we have $\lim_{h\rightarrow 0}A_h=0$. The second part of {\bf SC3} for $f$ follows from these properties and the general Leibniz rule.$\square$

\section{Appendix}\label{Appendix}

\begin{prop}\label{expo}
Assume that Assumptions {\bf E21-E22} hold true. Additionally, suppose that there exists $\kappa>0$ such that 
$$\int\exp(\kappa\vert y\vert)\left\vert f_{\varepsilon}^{(s)}(y)\right\vert dy<\infty,\quad s=0,\ldots,k.$$
There exists $\kappa'\in (0,\kappa)$ such that, setting $q=0$, $d_0=1$ and $\phi(y)=\exp(\kappa'\vert y\vert)$
, the conclusions of Corollary \ref{pMarkov} hold true.
\end{prop}

\paragraph{Proof of Proposition \ref{expo}}
We use Corollary \ref{pMarkov}.
We set $\sigma_{+}=\sup_{u\in [0,1]}\sigma(u)$ and $\overline{\beta}=\sup_{u\in [0,1]}\beta_{i,u}$.
We first fix $\kappa'\in (0,\kappa)$ small enough such that $\kappa'<\kappa\max\left((1-\overline{\beta})/\sigma_+,1\right)$.
We will not check the small set condition in {\bf SCp1}, the proof being similar to that of Proposition \ref{autoreg}. To check the drift condition, we use {\bf E21} and convexity of the exponential function.
\begin{eqnarray*}
R_u\phi({\bf x})&\leq& \E\left[\exp\left(\kappa'\sum_{i=1}^p\beta_{i,u}\vert x_i\vert+\kappa'\sigma_{+}\vert\varepsilon_1\vert\right)\right]\\
&\leq& \sum_{i=1}^p\beta_{i,u}\exp\left(\kappa' \vert x_i\vert\right)+\E\left[\exp\left(\frac{\sigma_+\kappa'}{1-\overline{\beta}}\vert\varepsilon_1\vert\right)\right].
\end{eqnarray*}
This shows the drift condition.
Assumption {\bf SCp2} is automatically satisfied and it remains to check {\bf SCp3}.
To this end, we will use the expression (\ref{shape}) and in particular the following bound which can be obtained using {\bf E22},
\begin{equation}\label{goodbound}
\left\vert\partial_1^{(\ell)}f(u,{\bf x},y)\right\vert\leq C\left[f_{\varepsilon}\left(\frac{y-m(u,{\bf x})}{\sigma(u)}\right)+\sum_{s=0}^{\ell}\left\vert f^{(s)}_{\varepsilon}\left(\frac{y-m(u,{\bf x})}{\sigma(u)}\right)\right\vert\right]\cdot\left[\vert y-m(u,{\bf x})\vert^s+(1+\vert x\vert^{q'})^s\right],
\end{equation}
for some constant $C>0$.
As for checking the drift condition, one can show that 
$$\int \exp(\kappa'\vert y\vert)f_{\varepsilon}\left(\frac{y-m(u,{\bf x})}{\sigma(u)}\right)dy\leq C\sum_{i=1}^p\phi(x_i),$$
for another positive constant $C$ (this constant can change from line to line).
Moreover, using our assumptions, we have two following inequalities,
\begin{eqnarray*} 
&&\int \exp\left(\kappa'\left\vert m(u,{\bf x})+\sigma(u)\varepsilon_1\right\vert\right)\left\vert f^{(s)}_{\varepsilon}(z)\right\vert\cdot\vert z\vert^s dz\\
&\leq& \sum_{i=1}^p\beta_{i,u}\exp(\kappa'\vert x_i\vert)\int \vert z\vert^s\cdot\left\vert f_{\varepsilon}^{(s)}(z)\right\vert dz+\int \exp\left(\frac{\sigma_+\kappa'\vert z\vert}{1-\overline{\beta}}\right)\left\vert f^{(s)}_{\varepsilon}(z)\right\vert\cdot\vert z\vert^sdz 
\end{eqnarray*}
and 
\begin{eqnarray*} 
&&(1+\vert {\bf x}\vert^{q'})^s\int \exp\left(\kappa'\left\vert m(u,{\bf x})+\sigma(u)\varepsilon_1\right\vert\right)\left\vert f^{(s)}_{\varepsilon}(z)\right\vert dz\\
&\leq& (1+\vert {\bf x}\vert^{q'})^s\left[\sum_{i=1}^p\beta_{i,u}\exp(\kappa'\vert x_i\vert)\int \left\vert f_{\varepsilon}^{(s)}(z)\right\vert dz+\int \exp\left(\frac{\sigma_+\kappa'\vert z\vert}{1-\overline{\beta}}\right)\left\vert f_{\varepsilon}^{(s)}(z)\right\vert dz\right].
\end{eqnarray*}
Using (\ref{goodbound}) and the condition on $\kappa'$, we get the first part of {\bf SCp3}.
To get the second part, it is sufficient to show that for $\ell=0,\ldots,k$,
$$\lim_{M\rightarrow \infty}\sup_{u\in [0,1]}\int_{\{\vert y\vert\geq M\}}\phi(y)\left\vert\partial_1^{(\ell)}
f(u,{\bf x},y)\right\vert dy=0.$$
This follows from some bounds that are similar to the previous ones and the condition on $\kappa'$. Details are omitted.$\square$

\begin{prop}\label{compar}
Assume that for each $u\in [0,1]$, there exist $C_u>0$ and $\kappa_u\in (0,1)$ such that for all $x\in G$,
$$\Vert \delta_xP_u^n-\mu_u\Vert_V\leq C_u\kappa_u^n.$$
Assume furthermore that for all $u\in [0,1]$, $\lim_{h\rightarrow 0}\Vert T_{u+h}-T_u\Vert_{V,V}=0$ and $\Vert T_0\Vert_{V,V}<\infty$,
where $T_u:\mathcal{M}_V(G)\rightarrow \mathcal{M}_V(G)$ is defined by $T_u\mu=\mu P_u$, $\mu\in\mathcal{M}_V(G)$.
Then the family of Markov kernel $\left\{P_u:u\in [0,1]\right\}$ is simultaneously $V-$uniformly ergodic.
\end{prop}

\paragraph{Proof of Proposition \ref{compar}}
First, we note that under our continuity assumption, we have $g:=\sup_{u\in [0,1]}\Vert T_u\Vert_{V,V}<\infty$.
Let $u\in [0,1]$. We have 
$$\Delta_V(P_u^n)\leq \sup_{x\in G}\frac{\Vert \delta_xP_u^n-\mu_u\Vert_V}{V(x)}\leq C_u\kappa_u^n.$$
The proof of the first inequality is given in the proof of Proposition \ref{central}.
There then exists an integer $n_u\geq 1$ such that $\Delta_V(P_u^{n_u})=\Vert T_u^{n_u}\Vert_{0,V,V}<1$.
By continuity of the application $u\mapsto T_u$, there exists a neighborhood $\mathcal{O}_u$ of $u$ such that
$\kappa_u:=\sup_{v\in \mathcal{O}_u}\Delta_V\left(P_v^{n_u}\right)<1$. 
Next, we show that $v\mapsto \pi_v$ is continuous at point $u$. 
If $v\in\mathcal{O}(u)$, we use the decomposition $\pi_v-\pi_u=(I-T_v)^{-1}(T_v-T_u)\mu_u$. This decomposition is given in the proof of Theorem \ref{general} (3.). One can note that
$$\kappa_u<1, g<\infty \Rightarrow \sup_{v\in\mathcal{O}_u}\Vert (I-T_u)^{-1}\Vert_{0,V,V}<\infty,$$   
as we showed in the proof of Theorem \ref{general} (with $[0,1]$ replaced by $\mathcal{O}_u$).
We then get
$$\Vert \nu_v-\nu_u\Vert_V\leq \sup_{v\in\mathcal{O}_u}\Vert (I-T_v)^{-1}\Vert_{0,V,V}\Vert T_v-T_u\Vert_{V,V}\Vert \mu_u\Vert_V.$$
From our continuity assumption on $u\mapsto T_u$, we deduce the continuity of $v\mapsto \pi_v$ at point $u$.
Since our result is valid for any $u\in [0,1]$, we deduce that $u\mapsto \mu_u$ is continuous and then $\sup_{u\in [0,1]}\Vert\mu_u\Vert_V<\infty$. Next, we will prove that for any $u\in [0,1]$, the family $\left\{P_v; v\in\mathcal{O}_u\right\}$ satisfies a simultaneously $V-$uniform ergodicity condition. Let $v\in [0,1]$, $n$ a positive integer, $s_u$ the integer part of the ration $n/n_u$ and $r_u=n-s_u n_u$. 
We have 
\begin{eqnarray*}
\Vert \delta_x P_v-\mu_v\Vert_V &=& \Vert \delta_x P_v^{r_u}P_v^{n_u s_u}-\mu_v P_v^{r_u}P_v^{n_u s_u}\Vert_V\\
&\leq& \Vert \delta_x P_v^{r_u}-\mu_v P_v^{r_u}\Vert_V\cdot \kappa_u^{s_u}\\
&\leq& \Vert \delta_x-\mu_v\Vert_V\cdot \Vert T_v\Vert^{r_u}_{V,V}\cdot \kappa_u^{s_u}.
\end{eqnarray*}
Setting $\overline{C}_u=\kappa_u^{-1}\Vert T_u\Vert_{V,V}\left(1+\sup_{v\in [0,1]}\Vert \mu_v\Vert_V\right)$
and $\overline{\kappa}_u=\kappa_u^{1/n_u}$, we have 
$$\Vert \delta_x P_v^n-\mu_v\Vert_V\leq \overline{C}_u V(x)\overline{\kappa}_u,$$
which shows this local simultaneous $V-$uniform ergodicity property.  
Its extension to the entire interval $[0,1]$ follows from a compactness argument, as in the proof of 
Proposition \ref{central}. $\square$

\paragraph{Acknowledgments.}  The author would like to thank Loïc Hervé and James Ledoux for some clarifications about the perturbation theory of Markov operators as well as four anonymous referees 
for many comments and suggestions that helped to improve considerably a first version of this paper.

\bibliographystyle{plainnat}
\bibliography{bibregul}

\end{document}